\documentclass{emsprocart}
\usepackage{mathrsfs,amssymb,amsfonts,amsmath}

\usepackage{graphics}
\usepackage[all]{xy}
\xyoption{curve}
\xyoption{import}
\xyoption{arc}
\xyoption{ps}
\UsePSspecials{dvips}
\usepackage{amsmath}
\usepackage{graphicx}
\usepackage[all]{xy}
\xyoption{curve}
\xyoption{import}
\xyoption{arc}
\xyoption{ps}
\UsePSspecials{dvips}

\contact[carrillo@math.jussieu.fr]{Paulo Carrillo Rouse\\ Projet
  d'alg{\`e}bres d'op{\'e}rateurs\\ Universit{\'e} de Paris 7\\ 175, rue de
  Chevaleret\\ 
Paris, France}

\newtheorem{theorem}{Theorem}[section]

\newtheorem{corollary}[theorem]{Corollary}

\newtheorem{lemma}[theorem]{Lemma}

\newtheorem{proposition}[theorem]{Proposition}

\newtheorem{ejemplo}{Example}
\newtheorem{teorema}{Theorem}
\newtheorem*{teorema1}{Theorem}

\theoremstyle{definition}

\newtheorem{definition}[theorem]{Definition}

\newtheorem{remark}[theorem]{Remark}

\newcommand{\ci}{C^{\infty}}

\newcommand{\Dnc}{\mathscr{D}}
\newcommand{\F}{\mathscr{F}}
\newcommand{\gr}{\mathscr{G}}
\newcommand{\go}{\mathscr{G} ^{(0)}}
\newcommand{\hr}{\mathscr{H}}
\newcommand{\ho}{\mathscr{H} ^{(0)}}
\newcommand{\gd}{\mathscr{G}^{\mathbb{R}^2}}
\newcommand{\gt}{\mathscr{G} ^{T}}
\newcommand{\Nb}{\mathscr{N}}
\newcommand{\Kom}{\mathscr{K}}
\newcommand{\ops}{\mathscr{O}}
\newcommand{\sw}{\mathscr{S}}
\newcommand{\Uo}{\mathscr{U}}
\newcommand{\Vo}{\mathscr{V}}
\newcommand{\Rr}{\mathbb{R}}
\newcommand{\Nat}{\mathbb{N}}
\newcommand{\src}{\mathscr{S}_{c}}
\newcommand{\cg}{C_{c}^{\infty}(\gr)}
\newcommand{\cgo}{C_{c}^{\infty}(\go)}
\newcommand{\ct}{C_{c}^{\infty}(\gr^T)}

\newcommand{\ckt}{C_{c}^{k}(\gr \times [0,1])}
\newcommand{\ck}{C_{c}^{k}(\gr)}

\title[Compactly supported analytic indices for Lie groupoids]{Compactly supported analytic indices for Lie groupoids}

\author{Paulo Carrillo Rouse}

\begin{document}

\begin{abstract}
For any Lie groupoid we construct an analytic index morphism taking values in a modified $K-theory$ group which involves the convolution algebra of compactly supported smooth functions over the groupoid. The construction is performed by using the deformation algebra of smooth functions over the tangent groupoid constructed in \cite{Ca2}. This allows in particular to prove a more primitive version of the Connes-Skandalis Longitudinal index Theorem for foliations, that is, an index theorem taking values in a group which pairs with Cyclic cocycles. As other application, for $D$ a
$\gr-$PDO elliptic operator with associated index
$ind\, D\in K_0(\ci_c (\gr))$, we prove that the
pairing $$<ind\, D,\tau>,$$ with $\tau$ a bounded continuous cyclic cocycle, 
only depends on the principal symbol class $[\sigma(D)]\in
K^0(A^*\gr)$. The result is completely general for {\'E}tale groupoids. We discuss some potential applications to the Novikov's conjecture.
\end{abstract}



\maketitle

\section{Introduction}
Index theory has as a starting point the Atiyah-Singer index theorem \cite{AS}: Let $D$ be an elliptic operator over a compact smooth manifold $M$, in particular $D$ has finite dimensional Kernel and Cokernel and the Fredholm index can be defined as $$ind\, D :=dim\, Ker\, D -dim\, Coker\, D\in \mathbb{Z}.$$
The map $D\mapsto ind\, D$ is completely codified by a group morphism
$K^0(T^*M)\longrightarrow \mathbb{Z},$ 
called the analytic index of $M$. That is, if $Ell(M)$ denotes the set of elliptic pseudodifferential operators over $M$, then the following diagram is commutative:
\begin{equation}\label{elldiam}
\xymatrix{
Ell(M)\ar[r]^-{ind} \ar[d]_-{\sigma}& \mathbb{Z}\\
K^0(T^*M)\ar[ru]_-{ind_{a,M}} & ,
}
\end{equation}
where $Ell(M)\stackrel{\sigma}{\longrightarrow}K^0(T^*M)$ is the surjective map that associates the class of the principal symbol of the operator in $K^0(T^*M)$. This is a fundamental property because it allows to use the (cohomological) properties of  $K-$theory and gives stability to the index. The topological formula given by Atiyah-Singer for the analytic index is the following:
\begin{equation}\label{topch}
ind_{a,M}([\sigma_P])=\int_{T^*M}ch([\sigma_P])Td(M).
\end{equation}
This kind of formulas give very interesting invariants of the manifold \cite{Law,Mel,Sha}.

We discuss now the Lie groupoids case. This concept is central in non commutative geometry. Groupoids
generalize the concepts of spaces, groups and 
equivalence relations. In the late 70's, mainly with the work
of Alain Connes, it became clear that groupoids appeared naturally as
substitutes of singular spaces \cite{Coinc,Mac,Ren, Pat}. Furthermore, Connes showed that many
groupoids and algebras associated to them appeared as `non commutative
analogues` of smooth manifolds to which many tools of geometry such as
K-theory and Characteristic classes could be applied \cite{Coinc,Concg}. 
Lie groupoids became a very natural place where to perform pseudodifferential calculus and index theory, \cite{Coinc,MP,NWX}. 

The study of the indices in the groupoid case is, as we will see, more delicate than the classical case. There are new phenomena appearing.
If $\gr$ is a Lie groupoid, a $\gr$-pseudodifferential operator is a differentiable family (see \cite{MP,NWX}) of operators. Let $P$ be such an operator, the index of $P$, $ind\, P$, is an element of $K_0(\cg)$. 
There is a diagram similar to (\ref{elldiam}) for the Lie groupoid setting:
\begin{equation}\label{ellgdiam}
\xymatrix{
Ell(\gr) \ar[r]^-{ind}\ar[d]_{symb}& K_0(\cg) \\
K^0(A^*\gr)  \ar@{.>}[ur]_{\nexists} & .
}
\end{equation}
We would like, as above, that the map $D\mapsto ind\, D\in
K_0(\cg)$ factor through $K^0(A^*\gr)$: this group codifies the principal symbols of all the elliptic $\gr$-operators and carries the topological invariants. Unfortunately, the map $D \mapsto ind\, D \in K_0(\cg)$ does not factor in general through the principal symbol \cite{Concg}. 
Nevertheless, if we consider the $K$-theory morphism 
$K_0(\cg)\stackrel{j}{\longrightarrow} K_0(C^*_r(\gr))$ 
induced by the inclusion $\cg \subset C^*_r(\gr)$, then the composition
$$Ell(\gr)\stackrel{ind}{\longrightarrow}K_0(\cg)\stackrel{j}{\longrightarrow} 
K_0(C^*_r(\gr))$$ do factor through the principal symbol by a morphism $ind_a$, called the analytic index of $\gr$:
\begin{equation}\label{ellgdiamc}
\xymatrix{
Ell(\gr) \ar[r]^-{ind}\ar[d]_{symb}& K_0(\cg)\ar[d]^-j \\
K^0(A^*\gr) \ar[r]_-{ind_a} & K_0(C^*_r(\gr)).
}
\end{equation}

In fact, the $C^*$-analytic index $ind_a$ can be completely described using the Connes tangent groupoid, \cite{Concg,MP,MN,DLN,ANS,HS}. We briefly recall this fact.

The tangent groupoid associated to a Lie groupoid $\gr \rightrightarrows \go$ is a Lie groupoid
$$\gr^T:=A\gr \times \{ 0\} \bigsqcup \gr \times (0,1]
\rightrightarrows \go \times [0,1],$$
that is compatible with $A\gr$ and $\gr$. Its $C^*$algebra is a continuous field of $C^*-$algebras over  
$[0,1]$ whose fiber at zero is 
$C_r^*(A\gr)\cong C_0(A^*\gr)$ and 
$C_r^*(\gr)$ elsewhere, it is a $C^*$-algebraic deformation in the sense of \cite{Land}. 
In particular it gives a short exact sequence of $C^*-$algebras
\begin{equation}\label{csegt}
0 \rightarrow C_r^*(\gr \times (0,1]) \longrightarrow C_r^*(\gr^T)
\stackrel{ev_0}{\longrightarrow}
C_r^*(A\gr) \rightarrow 0,
\end{equation}
Now, thanks to the good properties of the $K-$theory for $C^*$-algebras, Monthubert and Pierrot show in \cite{MP} that
\begin{equation}
ind_a=(ev_1)_*\circ
(ev_0)_*^{-1},
\end{equation}
where $C_r^*(\gr^T)
\stackrel{ev_1}{\longrightarrow}
C_r^*(\gr)$ is the evaluation at $1$. 


If we try to adapt the above arguments for the index $ind \, D\in K_0(\cg)$, we find several difficulties:
For instance, the map $K_0(\ci_c(A\gr))\longrightarrow K_0(C^*_r(A\gr))=K^0(A^*\gr)$ is not an isomorphism, {\it i.e.,}
$K_0(\ci_c(A\gr))$ is not the group that we would like. On the other hand, 
the kernel of the evaluation map $\ci_c(\gr^T)\stackrel{ev_0}{\longrightarrow}\ci_c(A\gr)$ 
is not easy to describe, and we cannot conclude its $K_0$ group vanishes since the $K$-theory for general topological algebras is not necessarily homotopy invariant. In general, the morphism 
$K_0(\cg)\stackrel{j}{\rightarrow} K_0(C^*_r(\gr))$ is not an isomorphism. A very simple example of this situation is $\Rr \rightrightarrows \{*\}$ with the vectorial group structure (\cite{Concg}, p.142).
When the groupoid is proper we do not perceive this nuance: it is an
isomorphism. This is because $\cg$ is stable under holomorphic calculus 
(\cite{CMR,KarK}).  
In general, the $K$-theory groups as $K_0(\cg)$ do not satisfy in general some fundamental properties as homotopy invariance and Bott periodicity. 



Our first step is to consider a "good quotient" of $K_0(\cg)$ where we are going to be able to factorize $[ind\, D]\in K_0(\cg)/\sim$ through the principal symbol, using the Connes tangent groupoid.

Let us note, for each $k\in \mathbb{N}$, $K_{0}^{h,k}(\gr)$ the quotient group of $K_0(\ck)$ by the equivalence relation induced by $\tiny{K_0(\ckt)\overset{e_0}{\underset{e_1}{\rightrightarrows}}K_0(\ck)}$. Let $K_{0}^{F}(\gr)= \varprojlim_{k}K_{0}^{h,k}(\gr)$ be the projective limit relative to the inclusions $\ck \subset
C_{c}^{k-1}(\gr)$. 
The main result of this paper is the following.

\begin{teorema}\label{elteointro}
There exists a group morphism
\begin{equation}\nonumber
ind_{a}^{F}: K^0(A^*\gr)\rightarrow K_{0}^{F}(\gr)
\end{equation}
such that the following diagram is commutative
\begin{equation}
\xymatrix{
Ell(\gr) \ar[d]_{symb} \ar[r]^-{ind}&
K_0(\ci_c(\gr)) \ar[d] & \\ K_0(A^*\gr) \ar[d]_{id} \ar[r]^{ind_{a}^{F}}& K_{0}^{F}(\gr) \ar[d] & \\
K_0(A^*\gr) \ar[r]_-{ind_a} & K_0(C_{r}^{*}(\gr)) &
}
\end{equation}
\end{teorema}

The index of the above theorem is called {\it{``the compactly supported analytic index of $\gr$"}}. 
The construction of these indices is also based on the Connes tangent groupoid, as in case of $C^*$-algebras. In fact, in \cite{Ca2} we constructed an algebra of smooth functions over the tangent groupoid that will allow in this paper to perform the compactly supported indices as "deformations". We recall the main result of \cite{Ca2}.
  
\begin{teorema1}
There exists an intermediate algebra $\src (\gr^T)$ consisting of smooth functions over the tangent groupoid
$$\ct \subset \src (\gr^T) \subset C_r^*(\gr^T),$$ such that it is a field of algebras over $[0,1]$, whose fibers are
$$\sw (A\gr) \text{ at } t=0, \text{ and }$$
$$\cg \text{ for } t\neq 0.$$  
\end{teorema1}  

The compactly supported index fits in a commutative diagram of the following type:
\[
\xymatrix{
K_0(\src (\gr^T)) \ar[d]_-{e_{1}^{F}} \ar[r]^{e_0}
& K_0(\sw (A\gr))  \ar[ld]^-{ind_{a}^{F}}&  \\
 K_{0}^{F}(\gr) & & .
}
\]
\subsection*{Applications}
\subsubsection*{Longitudinal index theorem}
First, we slightly modify the indices defined above. 
We take $K_{0}^{B}(\gr)=\varinjlim_{m} K_{0}^{F}(\gr \times \Rr^{2m})$, the inductive limit induced by the Bott morphisms $K_{0}^{F}(\gr \times \Rr^{2m})\stackrel{Bott}{\longrightarrow}
K_{0}^{F}(\gr \times \Rr^{2(m+1)})$. Let $$ind_{a,\gr}^{B}:K^0(A^*\gr)\rightarrow K_{0}^{B}(\gr)$$ be the morphism given by the composition of $ind_{a,\gr}^{F}:K^0(A^*\gr)\rightarrow
 K_{0}^{F}(\gr)$ followed by $K_{0}^{F}(\gr)
 \stackrel{Bott}{\longrightarrow}K_{0}^{B}(\gr)$. This  morphism satisfies also theorem \ref{elteointro}, and so it is called {\it the Periodic compactly supported index} of $\gr$.

The periodic index can be calculated by topological methods in the case of foliations. That is, we have a longitudinal index theorem that reinforce the Connes-Skandalis one, \cite{CS}. Indeed, we will show that the equality between the analytic and the topological indices takes already place in $K_{0}^{B}(\gr)$.

\begin{teorema}\label{tilintro}
Let $(M,F)$ be foliated manifold with holonomy groupoid $\gr$. Then, the topological index ({\`a} la Connes-Skandalis) can be defined with values in $K_{0}^{B}(\gr)$ and it coincides with the Periodic compactly supported index of $\gr$:
$$ ind_{a,\gr}^{B}=ind_{t,\gr}^{B}.$$
\end{teorema}

The reinforced longitudinal index theorem allows us to define 
an assembly map in our setting:
\begin{equation}\label{muF}
\mu_F:K_{*,\tau}(B\gr)\rightarrow K_{0}^{B}(\gr),
\end{equation}
given by 
$\mu_F(\delta_D)=ind_{a}^{B}(\sigma_D)$. It fits in the following commutative diagram
\begin{equation}\label{mus}
 \xymatrix{
K_{*,\tau}(B\gr) \ar[r]^-{\mu} \ar[rd]_-{\mu_F} & K_0(C_r^*(\gr)) \\
& K_{0}^{B}(\gr) \ar[u]_i,
}
\end{equation}
where $\mu$ is the classical Baum-Connes map, \cite{BCH,Tu3}.
\subsubsection*{Pairings with cyclic cohomology}

The interest to keep track on the $\ci_c$-indices is because at this level we can make explicit calculations via the Chern-Weil-Connes theory. In fact there is a pairing 
\cite{Concdg,Concg,Karhc}
\begin{equation}\label{accouplement}
\langle \_ \, , \_ \rangle :K_0(\cg)\times HP^*(\cg)\rightarrow \mathbb{C}
\end{equation}
There are several known cocycles over $\cg$. An important problem in Noncommutative Geometry is to calculate the above pairing in order to obtain numerical invariants from the indices in  $K_0(\cg)$, \cite{Concg,CMnov,GorLottfg}.

Now, we can expect an easy (topological) calculation only if the map 
$D\mapsto \langle D \, , \tau \rangle$ 
($\tau \in HP^*(\cg)$ fix) factors through the symbol class of $D$, $[\sigma(D)]\in K^0(A^*\gr)$: we want 
to have a diagram of the following kind:
\[
\xymatrix{
Ell(\gr) \ar[r]^-{ind} \ar[d]_{symb.} & K_0(\cg) \ar[rr]^-{\langle \_ , \tau \rangle} & & \mathbb{C} \\
K^0(A^*\gr) \ar@{.>}[urrr]_-{\tau} & &  &.
}
\]

The next step in this work will consist of solving this factorization problem, using the compactly supported indices. In fact, the Chern-Connes theory applies naturally to algebras as $\ck$. Moreover, the pairing with Periodic cyclic cohomology preserves the relation $\sim_h$ and is compatible with Bott morphism. The result is the following:

\begin{teorema}[Factorization theorem]
Let $\tau $ be a bounded continuous cocycle. Then $\tau$ defines a morphism 
$\Psi_{\tau}:K^0(A^*\gr)\rightarrow \mathbb{C}$ such that the following diagram is commutative
\[
\xymatrix{
Ell(\gr) \ar[r]^-{ind} \ar[d]_-{symb.} & K_0(\cg)  \ar[rr]^-{\langle \_ , \tau \rangle} & & \mathbb{C} \\
K^0(A^*\gr) \ar[urrr]_-{\Psi_{\tau}} &  &.
}
\]
\end{teorema}

The restriction of taking bounded continuous cyclic cocycles in the last theorem is not at all restrictive. In fact, all the geometrical cocycles are of this kind (Group cocycles, The transverse fundamental class, Godbillon-Vey and all the Gelfand-Fuchs cocycles for instance). Moreover, for the case of {\'e}tale groupoids, the  explicit calculations of the Periodic cohomologies spaces developed in \cite{BN,Cra} allow us to conclude that the above result is completely general in this setting (see Theorem \ref{extfinteo}). 
In the last part of this work we will discuss how to recover some known index formulas for foliations.

\noindent
{\bf Acknowledgments} 
The present work is part of my PHD thesis at the University of Paris 7, \cite{Ca3}. I would like to thank my advisor, Professor Georges Skandalis, for proposing me this subject and for his guide and support.

\section{Lie groupoids}

Let us recall what a groupoid is:

\begin{definition}
A $\it{groupoid}$ consists of the following data:
two sets $\gr$ and $\go$, and maps
\begin{itemize}
\item[$\cdot$]$s,r:\gr \rightarrow \go$ 
called the source and target map respectively,
\item[$\cdot$]$m:\gr^{(2)}\rightarrow \gr$ called the product map 
(where $\gr^{(2)}=\{ (\gamma,\eta)\in \gr \times \gr : s(\gamma)=r(\eta)\}$),
\end{itemize}
such that there exist two maps, $u:\go \rightarrow \gr$ (the unit map) and 
$i:\gr \rightarrow \gr$ (the inverse map),
such that, if we note $m(\gamma,\eta)=\gamma \cdot \eta$, $u(x)=x$ and 
$i(\gamma)=\gamma^{-1}$, we have 
\begin{itemize}
\item[1.]$r(\gamma \cdot \eta) =r(\gamma)$ and $s(\gamma \cdot \eta) =s(\eta)$.
\item[2.]$\gamma \cdot (\eta \cdot \delta)=(\gamma \cdot \eta )\cdot \delta$, 
$\forall \gamma,\eta,\delta \in \gr$ when this is possible.
\item[3.]$\gamma \cdot x = \gamma$ and $x\cdot \eta =\eta$, $\forall
  \gamma,\eta \in \gr$ with $s(\gamma)=x$ and $r(\eta)=x$.
\item[4.]$\gamma \cdot \gamma^{-1} =u(r(\gamma))$ and 
$\gamma^{-1} \cdot \gamma =u(s(\gamma))$, $\forall \gamma \in \gr$.
\end{itemize}
Generally, we denote a groupoid by $\gr \rightrightarrows \go $.
\end{definition}

Along this paper we will only deal with Lie groupoids, that is, 
a groupoid in which $\gr$ and $\go$ are smooth manifolds (possibly with boundary), and $s,r,m,u$ are smooth maps (with s and r submersions, see \cite{Mac,Pat}). For $A,B$ subsets of $\go$ we use the notation
$\gr_{A}^{B}$ for the subset $\{ \gamma \in \gr : s(\gamma) \in A,\, 
r(\gamma)\in B\}$.

We recall how to define an algebra structure in $\cg$ using
smooth Haar systems.
 
\begin{definition}
A $\it{smooth\, Haar\, system}$ over a Lie groupoid is a family of
measures $\mu_x$ in $\gr_x$ for each $x\in \go$ such that,

\begin{itemize}
\item for $\eta \in \gr_{x}^{y}$ we have the following compatibility
  condition:
$$\int_{\gr_x}f(\gamma)d\mu_x(\gamma)
=\int_{\gr_y}f(\gamma \circ \eta)d\mu_y(\gamma)$$
\item for each $f\in \cg$ the map
$$x\mapsto \int_{\gr_x}f(\gamma)d\mu_x(\gamma) $$ belongs to $\cgo$
\end{itemize}

\end{definition}

A Lie groupoid always posses a smooth Haar system. In fact, if we
fix a smooth (positive) section of the 1-density bundle associated to
the Lie algebroid we obtain a smooth Haar system 
in a canonical way. 
We suppose for the rest of the
paper a given smooth Haar system given by 1-densities (for complete
details see \cite{Pat}). 
We can now define a convolution
product on $\cg$: Let $f,g\in \cg$, we set

$$(f*g)(\gamma)
=\int_{\gr_{s(\gamma)}}
f(\gamma \cdot \eta^{-1})g(\eta)d\mu_{s(\gamma)}(\eta)$$

This gives a well defined associative product. 
\begin{remark}
There
is a way to define the convolution algebra
using half densities (see Connes book \cite{Concg}).
\end{remark}
As we mentioned in the introduction, we are going to consider some elements
in the $K$-theory group $K_0(\cg)$. We recall how these elements are
usually defined (See \cite{NWX} for complete details): First we recall a few facts about 
$\gr$-Pseudodifferential calculus:

A $\gr$-$\it{Pseudodifferential}$ $\it{operator}$ is a family of
pseudodifferential operators $\{ P_x\}_{x\in \go} $ acting in
$\ci_c(\gr_x)$ such that if $\gamma \in \gr $ and
$$U_{\gamma}:\ci_c(\gr_{s(\gamma)}) \rightarrow \ci_c(\gr_{r(\gamma)}) $$
the induced operator, then we have the following compatibility condition 
$$ P_{r(\gamma)} \circ U_{\gamma}= U_{\gamma} \circ P_{s(\gamma)}$$
There is also a differentiability condition with respect to $x$ that can
be found in \cite{NWX}.

For $P$ a $\gr$-Pseudodifferential elliptic operator there is a parametrix, $\it{i.e.}$, a $\gr$-Pseudodifferential
 operator $Q$ such that $PQ-1$ and $QP-1$ belong to $\cg$
(where we are identifying $\cg$ with $\Psi^{-\infty}(\gr)$ as in
\cite{NWX}). In other words, $P$ defines a quasi-isomorphism in 
$(\Psi^{+\infty},\cg)$ and thus an element in $K_0(\cg)$ that we call
the index $ind(P)$. Similarly to the classical case, a $\gr$-PDO
operator has a principal symbol that defines an element in the
$K$-theory group $K^0(A^*\gr)$. 

\section{Deformation to the normal cone}

The tangent groupoid is a particular case of a geometric construction that we describe here.

Let $M$ be a $\ci$ manifold and $X\subset M$ be a $\ci$ submanifold. We denote
by $\Nb_{X}^{M}$ the normal bundle to $X$ in $M$, $\it{i.e.}$, 
$\Nb_{X}^{M}:= T_XM/TX$.

We define the following set
\begin{align}
\Dnc_{X}^{M}:= \Nb_{X}^{M} \times {0} \bigsqcup M \times \Rr^* 
\end{align} 
The purpose of this section is to recall how to define a $\ci$-structure in $\Dnc_{X}^{M}$. This is more or less classical, for example
it was extensively used in \cite{HS}.

Let us first consider the case where $M=\Rr^p\times \Rr^q$ 
and $X=\Rr^p \times \{ 0\}$ (where we
identify canonically $X=\Rr^p$). We denote by
$q=n-p$ and by $\Dnc_{p}^{n}$ for $\Dnc_{\Rr^p}^{\Rr^n}$ as above. In this case
we clearly have that $\Dnc_{p}^{n}=\Rr^p \times \Rr^q \times \Rr$ (as a
set). Consider the 
bijection  $\psi: \Rr^p \times \Rr^q \times \Rr \rightarrow
\Dnc_{p}^{n}$ given by 
\begin{equation}\label{psi}
\psi(x,\xi ,t) = \left\{ 
\begin{array}{cc}
(x,\xi ,0) &\mbox{ if } t=0 \\
(x,t\xi ,t) &\mbox{ if } t\neq0
\end{array}\right.
\end{equation}
which inverse is given explicitly by 
$$
\psi^{-1}(x,\xi ,t) = \left\{ 
\begin{array}{cc}
(x,\xi ,0) &\mbox{ if } t=0 \\
(x,\frac{1}{t}\xi ,t) &\mbox{ if } t\neq0
\end{array}\right.
$$
We can consider the $\ci$-structure on $\Dnc_{p}^{n}$
induced by this bijection.

We pass now to the general case. A local chart 
$(\Uo,\phi)$ in $M$ is said to be a $X$-slice if 
\begin{itemize}
\item[1)]$\phi : \Uo \stackrel{\cong}{\rightarrow} U \subset \Rr^p\times \Rr^q$
\item[2)]If $\Uo \cap X =\Vo$, $\Vo=\phi^{-1}( U \cap \Rr^p \times \{ 0\}
  )$ (we note $V=U \cap \Rr^p \times \{ 0\}$)
\end{itemize}
With this notation, $\Dnc_{V}^{U}\subset \Dnc_{p}^{n}$ as an
open subset. We may define a function 
\begin{equation}\label{phi}
\tilde{\phi}:\Dnc_{\Vo}^{\Uo} \rightarrow \Dnc_{V}^{U} 
\end{equation}
in the following way: For $x\in \Vo$ we have $\phi (x)\in \Rr^p
\times \{0\}$. If we write 
$\phi(x)=(\phi_1(x),0)$, then 
$$ \phi_1 :\Vo \rightarrow V \subset \Rr^p$$ 
is a diffeomorphism. We set 
$\tilde{\phi}(v,\xi ,0)= (\phi_1 (v),d_N\phi_v (\xi ),0)$ and 
$\tilde{\phi}(u,t)= (\phi (u),t)$ 
for $t\neq 0$. Here 
$d_N\phi_v: N_v \rightarrow \Rr^q$ is the normal component of the
 derivate $d\phi_v$ for $v\in \Vo$. It is clear that $\tilde{\phi}$ is
 also a  bijection (in particular it induces a $C^{\infty}$ structure on $\Dnc_{\Vo}^{\Uo}$). 
Now, let us consider an atlas 
$ \{ (\Uo_{\alpha},\phi_{\alpha}) \}_{\alpha \in \Delta}$ of $M$
 consisting of $X-$slices. Then the collection $ \{ (\Dnc_{\Vo_{\alpha}}^{\Uo_{\alpha}},\tilde{\phi_{\alpha})}
  \} _{\alpha \in \Delta }$ is a $\ci$-atlas of
  $\Dnc_{X}^{M}$ (proposition 3.1 in \cite{Ca2}).

\begin{definition}[Deformation to the normal cone]
Let $X\subset M$ be as above. The set
$\Dnc_{X}^{M}$ equipped with the  $C^{\infty}$ structure
induced by the atlas described in the last proposition is called
$\it{"The\, deformation\, to\, normal\, cone\, associated\, to\,}$   
$X\subset M$". 
\end{definition}

\begin{remark}
Following the same steps, we can define a deformation to the normal
cone associated to an injective immersion $X\hookrightarrow M$.
\end{remark}

One important feature about this construction is that it is in
some sense functorial. More explicitly, let $(M,X)$ 
and $(M',X')$ be $\ci$-pairs as above and let
 $F:(M,X)\rightarrow (M',X')$
be a pair morphism, i.e., a $\ci$ map   
$F:M\rightarrow M'$, with $F(X)\subset X'$. We define 
$ \Dnc(F): \Dnc_{X}^{M} \rightarrow \Dnc_{X'}^{M'} $ by the following formulas:\\

$\Dnc(F) (x,\xi ,0)= (F(x),d_NF_x (\xi),0)$ and\\

$\Dnc(F) (m ,t)= (F(m),t)$ for $t\neq 0$,
\noindent
where $d_NF_x$ is by definition the map
\[  (\Nb_{X}^{M})_x 
\stackrel{d_NF_x}{\longrightarrow}  (\Nb_{X'}^{M'})_{F(x)} \]
induced by $ T_xM 
\stackrel{dF_x}{\longrightarrow}  T_{F(x)}M'$.

Then $\Dnc(F):\Dnc_{X}^{M} \rightarrow \Dnc_{X'}^{M'}$ is $\ci$-map (proposition 3.4 in \cite{Ca2}).

\subsection{The tangent groupoid}


\begin{definition}[Tangent groupoid]
Let $\gr \rightrightarrows \go $ be a Lie groupoid. $\it{The\, tangent\,
groupoid}$ associated to $\gr$ is the groupoid that has $\Dnc_{\go}^{\gr}$ as the set of arrows and  $\go \times \Rr$ as the units, with:
\begin{itemize}
\item[$\cdot$] $s^T(x,\eta ,0) =(x,0)$ and $r^T(x,\eta ,0) =(x,0)$ at $t=0$.
\item[$\cdot$] $s^T(\gamma,t) =(s(\gamma),t)$ and $r^T(\gamma,t)
  =(r(\gamma),t)$ at $t\neq0$.
\item[$\cdot$] The product is given by
  $m^T((x,\eta,0),(x,\xi,0))=(x,\eta +\xi ,0)$ et \linebreak $m^T((\gamma,t), 
  (\beta ,t))= (m(\gamma,\beta) , t)$ if $t\neq 0 $ and 
if $r(\beta)=s(\gamma)$.
\item[$\cdot$] The unit map $u^T:\go \rightarrow \gr^T$ is given by
 $u^T(x,0)=(x,0)$ and $u^T(x,t)=(u(x),t)$ for $t\neq 0$.
\end{itemize}
We denote $\gr^{T}:= \Dnc_{\go}^{\gr}$ and $A\gr:=\Nb_{\go}^{\gr}$.
\end{definition} 

As we have seen above $\gr^{T}$ can be considered as a $\ci$ manifold with
border. As a consequence of the functoriality of the Deformation to the normal cone,
one can show that the tangent groupoid is in fact a Lie
groupoid. Indeed, it is easy to check that if we identify in a
canonical way $\Dnc_{\go}^{\gr^{(2)}}$ with $(\gr^T)^{(2)}$, then 
$$ m^T=\Dnc(m),\, s^T=\Dnc(s), \,  r^T=\Dnc(r),\,  u^T=\Dnc(u)$$
where we are considering the following pair morphisms:
\begin{align}  
m:((\gr)^{(2)},\go)\rightarrow (\gr,\go ), \nonumber
\\
s,r:(\gr ,\go) \rightarrow (\go,\go),\nonumber 
\\
u:(\go,\go)\rightarrow (\gr,\go ).\nonumber
\end{align}

\begin{remark}\label{haartg}
Finally, let $\{ \mu_x\}$ be a smooth Haar system on $\gr$, {\it i.e.}, a choice of $\gr$-invariant Lebesgue measures. In particular we have an associated smooth Haar system on $A\gr$ (groupoid given by the vector bundle structure), which we note again by $\{ \mu_x\}$. Then the following family $\{\mu_{(x,t)}\}$ is a smooth Haar system for the tangent groupoid of $\gr$ (details may be found in \cite{Pat}):
\begin{itemize}
\item $\mu_{(x,0)}:=\mu_x$ at $(\gr^T)_{(x,0)}=A_x\gr$ and 
\item $\mu_{(x,t)}:=t^{-q}\cdot \mu_x$  at $(\gr^T)_{(x,t)}=\gr_x$ for
  $t\neq 0$, where $q=dim\, \gr_x$.
\end{itemize}
In this article, we are only going to consider this Haar systems for the tangent groupoids.
\end{remark}

\section{A Schwartz type algebra for the Tangent groupoid}

In this section we will recall how to construct the deformation algebra mentioned at the introduction.
For complete details, we refer the reader to \cite{Ca2}.

The Schwartz algebra for the Tangent groupoid will be a particular
case of a construction associated to any deformation to the normal
cone. 

\begin{definition}\label{ladef}
Let $p, q\in \Nat$ and $U \subset \Rr^p
  \times \Rr^q$ 
an open subset, and let $V=U\cap (\Rr^p \times \{ 0\})$.
\begin{itemize}
\item[(1)]Let $K\subset U \times \Rr$ be a compact
  subset. We say that $K$ is a conic compact subset of $U \times \Rr$
relative to $V$ if
\[ K_0=K\cap (U \times \{ 0\}) \subset V\]

\item[(2)]Let $\Omega_{V}^{U}=\{(x,\xi,t)\in 
\Rr^p \times \Rr^q \times \Rr: (x,t\cdot \xi)\in U \},$ 
which is an open subset of $\Rr^p \times \Rr^q \times \Rr$ and thus
a $\ci$ manifold. 
Let $g \in \ci (\Omega_{V}^{U})$. We say that
   $g$ has compact conic support, if there exists a conic
  compact $K$
 of $U \times \Rr$ relative to $V$ such that if 
$(x, t\xi ,t) \notin K$ then $g(x, \xi ,t)=0$.

\item[(3)]We denote by $\src (\Omega_{V}^{U})$ 
the set of functions
$g\in \ci (\Omega_{V}^{U})$ 
that have compact conic support and that satisfy the following condition:

\begin{itemize}
\item[$(s_1$)]$\forall$ $k,m\in \Nat$, $l\in \Nat^p$
and $\alpha \in \Nat^q$ it exists $C_{(k,m,l,\alpha)} >0$ such that
\[ (1+\| \xi \|^2)^k \| \partial_{x}^{l}\partial_{\xi}^{\alpha}
\partial_{t}^{m}g(x,\xi ,t) \| \leq C_{(k,m,l,\alpha)}   \]
\end{itemize}

\end{itemize}

\end{definition}

Now, the spaces $\src (\Omega_{V}^{U})$ are invariant under
diffeomorphisms. More precisely: Let $F:U\rightarrow U'$ be a $\ci$-diffeomorphism such that $F(V)=V'$; let 
$\tilde{F}:\Omega_{V}^{U}\rightarrow \Omega_{V'}^{U'}$ the induced map. Then, for every 
$g\in \src (\Omega_{V'}^{U'})$, we have that
$\tilde{g}:= g\circ \tilde{F} \in \src (\Omega_{V}^{U})$ (proposition 4.2 in \cite{Ca2}).

This compatibility result allows to give the following definition.

\begin{definition}\label{src}
Let $g \in \ci (\Dnc_{X}^{M}) $.
\begin{itemize}

\item[(a)]We say that $g$ has conic compact support $K$, if there exists a compact subset
 $K\subset M \times \Rr$ with $K_0:=K\cap (M\times \{ 0\}) \subset X$ (conic
 compact relative to $X$) such that if $t\neq 0$ and 
$(m,t) \notin K$ then $g(m,t)=0$.

\item[(b)]We say that $g$ is rapidly decaying at zero if for every
$(\Uo,\phi)$  $X$-slice chart
and for every $\chi \in \ci_c(\Uo \times \Rr)$, the map 
$g_{\chi}\in \ci(\Omega_{V}^{U})$ ($\Omega_{V}^{U}$ as in definition \ref{ladef}.)
given by
\[ g_{\chi}(x,\xi ,t)= (g\circ \varphi^{-1})(x,\xi ,t) 
\cdot (\chi \circ p \circ \varphi^{-1})(x,\xi ,t) \]
is in  $\src (\Omega_{V}^{U})$, where 
\begin{itemize}
\item[$\cdot$] $p$ is the projection is the deformation of the pair map $(M,X)\stackrel{Id}{\longrightarrow}
(M,M)$, {\it i.e.}, $p:\Dnc_{X}^{M} \rightarrow M\times \Rr$ is given by $(x,\xi,0)\mapsto (x,0)$, and 
$(m,t)\mapsto (m,t)$ for $t\neq 0$, and
\item[$\cdot$] $\varphi:=\tilde{\phi}^{-1}\circ \psi :\Omega_{V}^{U}\rightarrow \Dnc_{\Vo}^{\Uo}$, where $\psi$ and $\tilde{\phi}$ are defined at (\ref{psi}) and (\ref{phi}) above.
\end{itemize}
\end{itemize}

Finally, we denote by $\src (\Dnc_{X}^{M})$ the set of functions 
$g\in \ci(\Dnc_{X}^{M})$ that are rapidly decaying at zero 
with conic compact support.
\end{definition}

\begin{remark}\label{remsrc}
\noindent
\begin{itemize}
\item[{\it (a)}] Obviously $\ci_c(\Dnc_{X}^{M})$ is a subspace of $\src (\Dnc_{X}^{M})$.
\item[{\it (b)}] Let 
$\{ (\Uo_{\alpha},\phi_{\alpha}) \}_{\alpha \in \Delta}$ a family of
$X-$slices covering $X$. We have a decomposition of $\src(\Dnc_{X}^{M})$ as follows (see remark 4.5 in \cite{Ca2} and discussion below it):
\begin{align}\label{decomposicion}
\src (\Dnc_{X}^{M}) = \sum_{\alpha \in \Lambda} \src (\Dnc_{\Vo_{\alpha}}^{\Uo_{\alpha}}) 
+ \ci_c (M\times \Rr^*).
\end{align}
\end{itemize}
\end{remark}

The main theorem in \cite{Ca2} (Theorem 4.10) is the following

\begin{theorem}
The space $\src (\gr^T)$ is stable under convolution. In particular, we have the following inclusions of algebras
$$\ct \subset \src (\gr^T) \subset C_r^*(\gr^T)$$ 
Moreover, $\src (\gr^T)$ is a field of algebras over $\Rr$, whose fibers are
$$\sw (A\gr) \text{ at } t=0, \text{ and }$$
$$\cg \text{ for } t\neq 0.$$  
\end{theorem}  

In the statement of this theorem, $\sw (A\gr)$ denotes the Schwartz algebra over the Lie algebroid. Let us briefly recall the notion of Schwartz space associated to a
vector bundle: For a trivial bundle $X\times \Rr^q\rightarrow X$, $\sw(X\times \Rr^q):=\ci_c(X,\sw(\Rr^q))$ (see \cite{trev}). In general, $\sw(E)$ is defined using local charts. More precisely, a partition of the unity argument, allows to see that if we take a covering of $X$, $\{(\Vo_{\alpha},\tau_{\alpha})\}_{\alpha \in \Delta}$, consisting on trivializing charts, then we have a decomposition of the following kind:
\begin{equation}\label{des}
\sw(E)= \sum_{\alpha}\sw (\Vo_{\alpha}\times \Rr^q).
\end{equation}

The "Schwartz algebras" have in general the good $K-$theory
groups. As we said in te introduction, we are
interested in the group $K_0(A^*\gr)=K_0(C_0(A^*\gr))$. It is not 
enough to take the $K-$theory of $\ci_c(A\gr)$ (see for
example \cite{Concg}). We are going to see that $\sw (A^*\gr)$ has the wanted $K$-theory. In particular, our deformation algebra restricts at zero to the right algebra. The result is the following.

\begin{proposition}
Let $E\rightarrow X$ be a smooth vector bundle over a smooth manifold $X$.
\begin{itemize}
\item[{\it (i)}] The algebra $(\sw(E), *)$ with the convolution product is isomorphic to $(\sw(E^*),\cdot)$, where $\cdot$ denotes the punctual product.
\item[{\it (ii)}] The $\sw(E)$ is stable under holomorphic calculus on
$C_0(E^*)$. In particular, $K^0(E^*)\cong K_0(\sw(E))$.
\end{itemize}
\end{proposition}

\begin{proof}
\begin{itemize}
\item[{\it (i)}] Suppose first that $E=X\times \Rr^q$ is trivial. In this case $\sw(E)=\ci_c(X,\sw(\Rr^q))$.
Let $g\in \ci_c(X,\sw(\Rr^q))$ and $X\in \Rr^q$. We let
$$ \F(g) (x)(X)=\int_{\Rr^q}e^{-iX\cdot \eta}g(x,\eta)d\eta$$
which is the Fourier transform of $g(x)$ evaluated at $X$. It defines an element of
$\ci_c(X,\sw(\Rr^q))$. Now, since the product on $\sw(E)$ is given by
$$(f*g)(x)=f(x)*g(x),$$
we have, thanks to the continuity of the Fourier transform, that $\F$ is an isomorphism 
$(\ci_c(X,\sw(\Rr^q)),*)\cong (\ci_c(X,\sw(\Rr^q)),\cdot)$. 

For the general case, we look at the decomposition of $\sw(E)$ as in (\ref{des}) and we have the isomorphism given by Fourier
$$ \F: (\sw(E),* ) \rightarrow (\sw(E^*),\cdot).$$

\item[{\it (ii)}] Let $g\in \sw(E)$ and $f:\mathbb{C}\rightarrow \mathbb{C}$ an holomorphic map with 
$f(0)=0$. We have to verify that $f\circ g \in \sw(E)$: 
This condition is local, in particular we can suppose $E=V\times \Rr^q$ with $V\subset \Rr^p$ an open subset of $\Rr^p$. Let $\alpha\in \mathbb{N}^p \times \mathbb{N}^q$ and $K\subset V$ be a compact subset. We et $z=(x,\xi)\in V\times \Rr^q$  and $u=g(x,\xi)$. We can show by induction that the derivative
$\partial_{z}^{|\alpha|}(f\circ g)$ can be written in the following way:
$$\partial_{z}^{|\alpha|}(f\circ g)=\sum_{|\beta|\leq |\alpha|}\partial_{u}^{|\beta|}f(u)P_{\beta}(z),$$
where each $P_{\beta}(z)$ is a finite sum of products of the form
$$\partial_{z}^{\gamma_1}g(z)\cdots
\partial_{z}^{\gamma_r}g(z), $$ 
with $|\gamma_1|+...+|\gamma_r|=|\beta|$. Since $g\in \sw(E)$, $\partial_{u}^{|\beta|}f(u)$ is bounded for 
$x\in K$. The conclusion is immediate by using the Schwartz condition for $g$.
\end{itemize}
\end{proof}

From now on it will be important to restrict our functions on the tangent groupoid to the closed interval 
$[0,1]$. We keep the notation $\src (\Dnc_{X}^{M})$ for the restricted space. All the results above remain true. So for instance $\src (\gr^T)$ is an algebra which is a field of algebras over the closed interval $[0,1]$ with 0-fiber $\sw(A\gr)$ and $\ci_c(\gr)$ otherwise.
Before starting with the construction of the indices, we need to have an exact sequence analog as the one used in the construction of the $C^*$-analytic indices (exact sequence (\ref{csegt})). The first step in this direction is the following proposition.

\begin{proposition}
The evaluation at zero, $\src(\Dnc_{X}^{M})\stackrel{e_0}{\longrightarrow}\sw(\Nb_{X}^{M})$ is surjective.
\end{proposition}

\begin{proof}
Thanks to the decomposition (\ref{decomposicion}) discussed at the remark \ref{remsrc} above, it will be enough to prove that the evaluation map 
$\src(\Omega_{V}^{U})\stackrel{e_0}{\longrightarrow}\sw(V\times \Rr^q)$ is surjective. Where we are using the same notations as in definition \ref{ladef}.

Let $g\in \sw (V\times \Rr^q)$. We consider
$$
h(s) = \left\{
\begin{array}{cc}
e^{-\frac{1}{s}} &\mbox{ if } s>0 \\
0 &\mbox{ if } s\leq0
\end{array}\right.
$$
Let $K_0\subset V$ be the horizontal support of $g$. We can assume without lost of generality that 
$\{ (x,\xi)\in \Rr^p \times \Rr^q : x\in K_0 \text{ and } \| \xi \| \leq 1 \} \subset U$.
Let
\[ K=\{ (x,\xi,t)\in U\times [0,1]: x\in K_0,\| \xi\| \leq \sqrt{t} \},\]
then $K\subset U\times [0,1]$ is a compact subset and 
$K\cap U\times \{ 0\} =K_0$, $\it{i.e.}$, $K$ is a conic compact $U\times [0,1]$ relative to $V$.

Let $(x,\xi,t)\in \Omega_{V}^{U}$, and
$$
\tilde{g}(x,\xi,t)= \left\{
\begin{array}{cc}
g(x,\xi)\cdot h(1-t\| \xi \|^2) &\mbox{ if } x\in K_0 \\
0 &\mbox{ otherwise }
\end{array}\right.
$$
then we have that

\begin{itemize}
\item[$\cdot$]$\tilde{g}$ has compact conic support $K$
\item[$\cdot$]$\tilde{g} \in \ci (\Omega_{V}^{U})$ because $g$ and $h$ are
  $\ci$ and $g$ has horizontal compact support contained in $K_0$.
\item[$\cdot$]$\tilde{g}\in \src (\Omega_{V}^{U})$: Let 
$\tilde{h}(\xi,t)=h(1-t\| \xi \|^2)$. By induction, an elementary computation shows that, for $\alpha \in \Nat^q$, $l\in \Nat$, we have 
\begin{center}
$\partial_{\xi}^{\alpha}\partial_{t}^{l}\tilde{h}(\xi,t)=
\sum_{|\beta |\leq |\alpha|+l}a_{\beta}P_{\beta}(\xi,t) h^{|\beta|}(1-t\| \xi \|^2)$,
\end{center}
where $\beta \in \Nat^{q+1}$, $a_{\beta}$ is a constant (depending on $\alpha$ and $l$) and $P_{\beta}(\xi,t)$ is a finite sum of products
\begin{center}
$t^{\gamma_0}\cdot \xi_{1}^{\gamma_1}\cdots \xi_{q}^{\gamma_q}$,
$\gamma_i\in \Nat$.
\end{center}
Now, from the fact that $h$ and all its derivates are bounded, and using the Schwratz property for $g$, we conclude that
$\tilde{g}\in \src (\Omega_{V}^{U})$. Finally, $e_0(\tilde{g})=g$ by construction.
\end{itemize}
\end{proof}

We have a the short exact sequence of algebras:
\begin{align}\label{se}
0 \longrightarrow J \longrightarrow \src (\gr^T)
\stackrel{e_0}{\longrightarrow} \sw (A\gr) \longrightarrow 0,
\end{align}
where $J=Ker(e_0)$ by definition.

\section{Compactly supported analytic indices}

This section is devoted to the construction of the
indices announced in the introduction.

As a first step, we want to apply $K$-theory to the exact sequence (\ref{se}) above. But in principle there is no reason for obtaining an exact sequence of the same kind. We have the following proposition:

\begin{proposition}\label{eosur}
The morphism in $K-$theory, 
\[ K_0(\src (\gr^T))
\stackrel{(e_0)_*}{\longrightarrow} K_0(\sw (A\gr))\approx K^0(A^*\gr),  \]
induced from the evaluation at zero is surjective.
\end{proposition}

\begin{proof}
Let $[\sigma] \in K_0(\sw (A\gr))=K^0(A^*\gr)$. 
We know from the $\gr$-pseudodifferential calculus that $[\sigma]$ can be represented by a smooth homogeneus elliptic symbol. We can consider the symbol over $A^*\gr
\times [0,1]$ that coincides with $\sigma$ for all $t$, we note it 
by $\tilde{\sigma}$. Now, since
$A\gr^T=A\gr
\times [0,1] $, we can take $\tilde {P}=(P_t)_{t\in [0,1]}$ a 
$\gr^T$-elliptic pseudodifferential operator associated to $\sigma$, that is,
$\sigma_{\tilde{P}}=\tilde{\sigma}$. Let
$i:\ci_c (\gr^T) \rightarrow \src (\gr^T)$ be the inclusion (which is an algebra morphism), then 
$i_*(ind\, \tilde{P}) \in K_0(\src
(\gr^T))$ is such that 
$e_{0,*}(i_*(ind\, \tilde{P}))=\sigma$.
\end{proof}

By applying $K-$theory to the exact sequence (\ref{se}) we obtain
\begin{align}\label{kse}
K_0(J) \longrightarrow K_0(\src (\gr^T))
\stackrel{e_0}{\longrightarrow} K^0(A^*\gr) \longrightarrow 0,
\end{align}

Let us now define the group where our indices take values. 

\begin{definition}\label{khk}
Let $k\in \mathbb{N}$. We note by
\[ K_{0}^{h,k}(\gr) :=
Coeq
\left[\tiny{K_0(\ckt)\overset{e_0}{\underset{e_1}{\rightrightarrows}}K_0(\ck)}
\right] \]
the co-equalizer of the $K$-theory morphisms induced by the evaluations at zero and at one. In other words, 
$ K_{0}^{h,k}(\gr) $ is just the quotient of $K_0(\ck)$ by the image of $(e_1-e_0):K_0(\ckt)\rightarrow K_0(\ck)$. We note by $\pi : K_0(\ck)\rightarrow K_{0}^{h,k}(\gr)$ the quotient morphism.
We define {\it{the Bounded $K$-theory group of $\gr$}} as the projective limit
$$K_{0}^{F}(\gr)=\varprojlim_k K_{0}^{h,k}(\gr)$$
induced by the inclusions $...\hookrightarrow \ck \hookrightarrow C_{c}^{k+1}(\gr)\hookrightarrow ...$
\end{definition}

\begin{remark}\label{homalg}
In the above definition, the superscript $h$ in $K_{0}^{h,k}(\gr)$ 
makes reference to homotopy. Indeed, we can see the co-equalizer 
$K_{0}^{h,k}(\gr)$ otherwise:
for $x,y\in K_0(\ck) $, we say the they are homotopic, $x\sim_hy$,
if there exists $z\in K_0(\ckt) $ such that
\begin{itemize}
\item[($0$)] $e_0(z)=x$ and
\item[($1$)] $e_1(z)=y$.
\end{itemize}
This gives an equivalence relation in
$K_0(\ck)$ compatible with the group structure. The quotient 
$K_0(\ck)/\sim_h$ coincides with $K_{0}^{h,k}(\gr)$.
\end{remark}

Before stating our main theorem, let us see that there are natural morphisms 
$$K_0(\cg)\rightarrow K_{0}^{F}(\gr) \rightarrow K_0(C^*_r(\gr)).$$
The first is induced from the canonical morphisms $K_0(\cg) \rightarrow K_{0}^{h,k}(\gr)$ by using the universal property of the projective limits. 
\noindent
For the second one, we see that
\begin{equation}\label{coeqlim}
K_{0}^{F}(\gr) \cong
Coeq
\left[\tiny{\varprojlim_kK_0(\ckt)
\overset{\tilde{e_0}}{\underset{\tilde{e_1}}{\rightrightarrows}}
\varprojlim_kK_0(\ck)}
\right]
\end{equation}
That is, $K_{0}^{F}(\gr)$ is also a coequalizer. Now, the morphism $K_{0}^{F}(\gr) \rightarrow K_0(C^*_r(\gr))$ is induced from the canonical morphisms
$K_{0}^{h,k}(\gr) \rightarrow K_0(C^*_r(\gr))$ (the $K-$theory for $C^*$-algebras respects the homotopy relation).

Our main result is the following:

\begin{theorem}\label{teo}
\quad
\begin{enumerate}
\item There is a unique group morphism
$$ ind_{a}^{F}: K^0(A^*\gr)\rightarrow K_{0}^{F}(\gr) $$
that fits in the following commutative diagram
\begin{equation}\label{eldiagramaevaluado}
\xymatrix{
&K_0(\src(\gr^T))\ar[ld]_{e_0} \ar[rd]^{e_1^F}&\\
K^0(A^*\gr)\ar[rr]_{ind_{a}^{F}}&&K_{0}^{F}(\gr). 
}
\end{equation}
\item The morphism $ ind_{a}^{F}: K^0(A^*\gr)\rightarrow K_{0}^{F}(\gr) $ fits in the following commutative diagram
\begin{equation}\label{eldiagrama}
\xymatrix{
Ell(\gr) \ar[d]_{symb} \ar[r]^{ind}&
K_0(\ci_c(\gr)) \ar[d] & \\ K_0(A^*\gr) \ar[d]_{id} \ar[r]^{ind_{a}^{F}}& K_{0}^{F}(\gr) \ar[d] & \\
K_0(A^*\gr) \ar[r]_{ind_a}& K_0(C_{r}^{*}(\gr)) & 
.}
\end{equation}
\end{enumerate}
\end{theorem}

The proof of the theorem will require several lemmas:

\begin{lemma}\label{lambda}
Let $\lambda: \gr^T \rightarrow [0,1]$ be the projection and  $J=Ker(e_0)$, where $\src (\gr^T)
\stackrel{e_0}{\longrightarrow} \sw (A\gr)$, as in 
(\ref{se}). Then $J=\lambda \cdot \src (\gt)$.
\end{lemma}
\begin{proof}
The inclusion $\lambda \cdot \src (\gt) \subset J$ is obvious.
Let $f\in \src(\gr^T)$ with $e_0(f)=0$. As a $\ci$-map, $f\in \ci(\gr^T)$, we can look at its Taylor expansion with respect to $t$. We have then that there is $g\in \ci(\gr^T)$ such that $f=\lambda \cdot g$. 
Now, in the definition of $\src (\gr^T)$ we imposed a condition on the partial derivates $\partial_t$ 
with respect to $t$ (condition $(s_1)$ in definition \ref{ladef}). This condition implies $g\in \src(\gr^T)$. That is, $J\subset  \lambda \cdot \src (\gr^T)$.
\end{proof}


\begin{lemma}\label{fkdcn}
Let $(M,X)$ be a $\ci$ pair.
Let $f\in  \src (\Dnc_{X}^{M})$ and $k\in \mathbb{N}^*$. We define, for
$(m,t)\in M\times [0,1]$,
$$
f_k(m,t) = \left\{
\begin{array}{cc}
0 &\mbox{ if } t=0 \\
t^kf(m,t) &\mbox{ if } t\neq 0
\end{array}\right.
$$
Then $f_k\in C_{c}^{k}(M\times [0,1])$.
\end{lemma}

\begin{proof}
Since $C^k$ is a local property, we can assume 
$M=U\subset \Rr^p\times \Rr^q$ open and $X=V=U\cap (\Rr^p \times \{ 0
\})$. We have then to show that if $f\in \src (\Omega)$ (we remind
$\Omega =\{ (x,\xi,t)\in \Rr^p\times \Rr^q \times [0,1]: (x,t\xi)\in U
\}$) then
$$
f_k(x,\xi,t) = \left\{
\begin{array}{cc}
0 &\mbox{ if } t=0 \\
t^kf(x,\frac{\xi}{t},t) &\mbox{ if } t\neq 0
\end{array}\right.
$$
belongs to $C_{c}^{k}(U \times [0,1])$. Let us show this point.

First, the only problem is at $t=0$, out of zero $f_k$ is $\ci$. Let
$(x,\xi,0) \in \Omega$, we are going to see that $f_k$ is $C^k$ at this point. We can suppose $\xi \neq 0$ because otherwise the result is trivial.

Let $u_{\xi}(t)=\frac{\xi}{t}$, we have $u_{\xi} \in \ci ((0,1],\Rr^q)$
and $u_{\xi}$ satisfies
\begin{equation}\label{ulema}
lim_{t\rightarrow 0} \| u_{\xi} (t) \| =+\infty
\end{equation}

We recall  $\src (\Omega_{V}^{U})$ consists of compact conic supported maps
$g\in \ci (\Omega_{V}^{U})$
satisfying the following condition:

\begin{itemize}
\item[$(r_1$)]$\forall$ $n,m\in \mathbb{N}$, $l\in \mathbb{N}^p$
and $\alpha \in \mathbb{N}^q$ it exists $C_{(n,m,l,\alpha)} >0$ such that
\[ (1+\| \xi \|^2)^n \| \partial_{x}^{l}\partial_{\xi}^{\alpha}
\partial_{t}^{m}g(x,\xi ,t) \| \leq C_{(n,m,l,\alpha)}   \]
\end{itemize}

By hypothesis $f\in \src(\Omega)$, hence, thanks to equation (\ref{ulema}) above, we have for that
\begin{equation}\label{limu}
lim_{t\rightarrow 0} \| u_{\xi}(t) \|^r
\| \partial_{z}^{\alpha}f(x,u_{\xi}(t),t)  \| = 0, 
\end{equation}
for every $r\in \Nat$ and $\alpha \in \mathbb{N}^p \times \mathbb{N}^q \times
\mathbb{N}$, here $z=(x,\xi,t)$.

By definition $f_k(x,\xi,t)=t^k\cdot f(x,u_{\xi}(t),t)$ for $t\neq
0$ and zero otherwise, in particular $f_k$ is $\ci$ for $t\neq 0$. 
Now, out of zero, an induction argument shows that for $\alpha=(l,\gamma) 
\in \mathbb{N}^p \times (\mathbb{N}^q \times
\mathbb{N})$ we have that
\[
\partial_{z}^{\alpha}f_k(x,\xi,t)= t^{k-|\gamma|}\sum_{|\beta|\leq |\alpha|}a_{\beta}P_{\beta}(x,\xi,t),
\]
where $a_{\beta}$ are constants (depeding on $\alpha$) and $P_{\beta}(x,\xi,t)$ is a finite sum of products of the following type
\begin{center}
$t^{m_0}\cdot (u_{\xi_{1}}(t))^{m_1}\cdots (u_{\xi_{q}}(t))^{m_q}\cdot 
\partial_{z}^{\delta}f(x,u_{\xi}(t),t)$,
$m_i\in \Nat$,
\end{center}
where $u_{\xi_{i}}(t)=\frac{\xi_i}{t}$. Then, if $|\gamma|\leq k$ we can use limit (\ref{limu}) to obtain that
\[
lim_{t\rightarrow 0} \| \partial_{z}^{\alpha}f_k(x,\xi,t)  \| = 0.
\]
That is, $f_k$ is at least of class $C^k$ at zero.
Finally, the support of $f_k$ is contained in the conic support of $f$, which is compact by assumption.
\end{proof}

\begin{lemma}
Let $k\in \mathbb{N}$ and $q:=dim\, \gr_x$, we define (for $J=Ker(e_0)$, where $\src (\gr^T)
\stackrel{e_0}{\longrightarrow} \sw (A\gr)$)
\[ \varphi^{k} : J^{k+q} \rightarrow \ckt \]
in the following way:
$$
\varphi^k (f)(\gamma,t) = \left\{
\begin{array}{cc}
0 &\mbox{ if } t=0 \\
t^{-q}f(\gamma,t) &\mbox{ if } t\neq 0
\end{array}\right.
$$

Then, we have a well defined morphism of algebras
\[ \varphi^{k} : J^{k+q} \rightarrow \ckt \]
\end{lemma}

\begin{proof}
The fact that $\varphi_k$ is well defined is an immediate consequence of lemmas \ref{lambda} and \ref{fkdcn} above. Let us check that it is an algebra morphism. We recall that we are only considering Haar systems of the type described at \ref{haartg} above.

Let $f,g \in \src (\gr^T)$, hence
$$
\hbox{ \footnotesize $
\varphi^k(f*g)(\gamma,t)= \left\{
\begin{array}{cc}
0 &\mbox{ if } t=0 \\
\\
t^{-q} \displaystyle \int_{\gr_{s(\gamma)}}f(\gamma \circ \eta^{-1})g(\eta,t)t^{-q}
d\mu_{s(\gamma)}(\eta) &\mbox{ if } t\neq 0
\end{array}\right.
$}
$$
$$
\hbox{ \footnotesize $
=\left\{
\begin{array}{cc}
0 &\mbox{ if } t=0 \\
\\
 \displaystyle \int_{\gr_{s(\gamma)}}t^{-q}\cdot f(\gamma \circ \eta^{-1})t^{-q}\cdot g(\eta,t)
d\mu_{s(\gamma)}(\eta) &\mbox{ if } t\neq 0
\end{array}\right.
$}
$$
$$
=(\varphi^k(f)*\varphi_k(g))(\gamma,t).
$$
\end{proof}

\begin{remark}\label{losdiagramas}
Let $\varphi_k : K_0(J^{k+q}) \rightarrow
K_0(C_{c}^{k}(\gr \times [0,1]))$
the $K$-theory morphism induced by $\varphi^{k}$. By construction we have the two following properties:
\begin{itemize}
\item[(a)] $(e_0)_* \circ \varphi_{k}=0$ o{\`u} where $e_0$ denotes the evaluation
$e_0:C_{c}^{k}(\gr \times [0,1])\rightarrow C_{c}^{k}(\gr) $.
\item[(b)] The following diagram is commutative
\[
\xymatrix{
0 & & & \\
K_0(J) \ar[r] \ar[u] &
K_0(\src (\gr^T)) \ar[r]^{e_0} \ar[dr]^{e_{1}^{k}}
& K_0(\sw (A\gr)) \ar[r] & 0\\
K_0(J^{k+q}) \ar[u]^{j} \ar[r]_-{\varphi_{k}} &
K_0(C_{c}^{k}(\gr \times [0,1])) \ar[r]_-{e_1} & K_0(C_{c}^{k}(\gr)) \\
& & &
}
\]
where $K_0(J^{k+q})\stackrel{j}{\longrightarrow}K_0(J)$ is induced by the inclusion $J^{k+q}\subset J$. The fact that $j$ is surjective is immediate from the $K$-theory exact sequence 
$$K_0(J^N) \rightarrow K_0(J) \rightarrow K_0(J/J^N),$$
since $K_0(J/J^N)=0$ because $J/J^N$ is a nilpotent algebra (see for instance \cite{Ros}).
\end{itemize}
\end{remark}

\begin{lemma}
Let $w\in K_0(\src (\gr^T))$ with $w\in Ker(e_0)_*$. Then $e_{1}^{h,k}(w)=0$, where $e_{1}^{h,k}$ is the composition
$$K_0(\src (\gr^T))\stackrel{e_1}{\rightarrow} K_0(\ci_c(\gr))
\stackrel{\iota_k}{\rightarrow} K_0(C_{c}^{k}(\gr))
\stackrel{\pi}{\rightarrow}
K_{0}^{h,k}(\gr).$$
\end{lemma}

\begin{proof}
By the exact sequence in $K$-theory, (\ref{kse}), 
we can take $x\in K_0(J)$ such that $i_*(x)=w$. We can choose 
$y\in K_0(J^{k+q})$ with $j(y)=x$, because $K_0(J^{k+q})\stackrel{j}{\longrightarrow}K_0(J)$ is surjective (see discussion above). Now, the condition (b) in the precedent remark implies
$e_1(\varphi_{k}(y))=e_{1}^{k}(w)$, and since $e_0(\varphi_{k}(y))=0$ it follows that 
$e_{1}^{k}(w)=e_1(\varphi_{k}(y))\sim_h0$ in $K_0(C_{c}^{k}(\gr))$.
\end{proof}

We can now pass to the proof of the theorem:

\begin{proof}[Proof of theorem \ref{teo}.]
\quad
\begin{enumerate}
\item We have only to show the existence of the mentioned morphism since uniqueness will follow immediately from the surjectivity of the evaluation morphism $K_0(\src(\gr^T))\stackrel{e_0}{\longrightarrow}K_0(A^*\gr)$, (proposition \ref{eosur}). 

Let $k\in \mathbb{N}$. Let $\sigma \in K_0(A^*\gr)$ and take $w_{\sigma}\in K_0(\src (\gr^T))$ 
with $e_0(w_{\sigma})=\sigma$. We  put 
$$ind_{a}^{h,k}(\sigma):= e_{1}^{h,k}(w_{\sigma}).$$
From lemma above it follows that $ind_{a}^{h,k}(\sigma)$ does not depend on the choice of 
$w_{\sigma}$. Finally, the fact that $ind_{a}^{h,k}$ is a group morphism follows immediately from the fact that $e_0:K_0(\src(\gr^T))\rightarrow K_0(A^*\gr)$ and $e_1^k:K_0(\src(\gr^T))\rightarrow K_{0}^{h,k}(\gr)$ are group morphisms.

Now, the morphism
$$ind_{a}^{F}: K^0(A^*\gr)\rightarrow K_{0}^{F}(\gr)$$ is just the induced by all the indices $ind_{a}^{h,k}$ by the universal property of the projective limit. 
\item It is enough to show the commutativity of the diagram (\ref{eldiagrama}) for each $ind_{a}^{h,k}$:

First, let us consider the set of elliptic operators on the tangent groupoid, $Ell(\gr^T)$. From this set, we also have two evaluation maps
$$Ell(A\gr)\stackrel{e_0}{\longleftarrow} Ell(\gr^T)\stackrel{e_1}{\longrightarrow} Ell(\gr).$$
It is immediate that this evaluation maps commute with indices, {\it i.e.}, the following diagrams are commutative
\[
\xymatrix{
Ell(A\gr)\ar[d]_{ind}&Ell(\gr^T)\ar[l]_{e_0}\ar[r]^{e_1}\ar[d]^{ind}&Ell(\gr)\ar[d]^{ind}\\
K_0(\ci_c(A\gr))&K_0(\ci_c(\gr^T))\ar[l]^{e_0}\ar[r]_{e_1}&K_0(\cg).
}
\]
Now, we know from the existence of a asymptotic pseudodifferential calculus that $e_1$ is a surjective function (see \cite{Connes_Higson,MP,NWX}). In fact, 
the symbol map $Ell(\gr)\stackrel{\sigma}{\rightarrow}K^0(A^*\gr)$ can (alternately) defined as follows: Let 
$P\in Ell(\gr)$ and $\tilde{P}\in Ell(\gr^T)$ an $e_1$-lifting ($e_1(\tilde{P})=P$). We set 
$$\sigma(P):=j(ind(e_0(\tilde{P}))),$$
where $Ell(A\gr)\stackrel{ind}{\rightarrow}K_0(\ci_c(A\gr))\stackrel{j}{\rightarrow}K^0(A^*\gr)$. That is, the symbol map fits in the following commutative diagram
\[
\xymatrix{
Ell(\gr^T)\ar[d]^{e_0}\ar[rr]^{e_1}&&Ell(\gr)\ar[d]^{symb.}\\
Ell(A\gr)\ar[rd]_{ind} \ar[rr]_{ind_{a,A\gr}}&&K^0(A^*\gr)\\
&K_0(\ci_c(A\gr))\ar[ru]_{j}&
}
\]

Now, for proving the commutativity of
\[
\xymatrix{
Ell(\gr) \ar[dd]_-{symb} \ar[r]^-{ind}&
K_0(\ci_c(\gr)) \ar[dd]^-{\pi \circ \iota}  \\ &   \\
K_0(A^*\gr) \ar[r]_-{ind_{a}^{h,k}}&  K_{0}^{h,k}(\gr)
}
\]
we decompose it in commutative diagrams in the following way
\[
\hbox{ \scriptsize $
\xymatrix{
Ell(\gr) \ar[dddd]_{symb} \ar@/^3pc/[rrrr]^{ind} 
&Ell(\gr^T)\ar[l]_{e_1}\ar[d]_{e_0}\ar[rr]^{ind} &
 & K_0(\ci_c(\gr^T))\ar[ddll]^{e_0}\ar[ddddl]^j\ar[r]^{e_1} & 
 K_0(\ci_c(\gr)) \ar[ldddd]^{j} \ar[dddd]^{\pi \circ j}\\
 &Ell(A\gr) \ar[d]_{ind} &&&\\
 &K_0(\ci_c(A\gr))\ar[ldd]_{j}  &   
 & & \\
 &&&& \\
K^0(A^*\gr) \ar@/_3pc/[rrrr]_{ind_{a}^{h,k}}& & K_0(\src(\gr^T)) \ar[ll]^{e_{0}}
\ar[r]_{e_{1}^{k}}
& K_0(\ck) \ar[r]_{\pi} & K_{0}^{h,k}(\gr)
}$
}
\]

Now, for the commutativity of
\[
\xymatrix{
K^0(A^*\gr) \ar[d]_{id}\ar[r]^{ind_{a}^{h,k}}& K_{0}^{h,k}(\gr) \ar[d]^{\iota} & \\
K^0(A^*\gr) \ar[r]_{ind_a}& K_0(C_{r}^{*}(\gr)) &
}
\]
we proceed as above, {\it{i.e.}}, we decompose it as
\[
\xymatrix{
K^0(A^*\gr) \ar[dd]_{id} & K_0(\src(\gr^T)) \ar[l]^{e_0} \ar[d]^{\iota} \ar[r]^{e_{1}^{k}} &
K_0(\ck) \ar[dd]^{\iota} \ar[rrd]^{\pi}
 &  & \\
 & K_{0}(C_{r}^{*}(\gr^T)) \ar[rd]^{e_1} &  & & K_{0}^{h,k}(\gr) \ar[lld]^{\tilde{\iota}}\\
K^0(A^*\gr) \ar[rr]_{ind_a } \ar[ru]^{e_{0}^{-1}} & & K_0(C_{r}^{*}(\gr)),
&  &
}
\]

The fact that the index $ind_{a}^{F}$ 
also satisfies diagram (\ref{eldiagrama}) comes from the universal properties of projective limits and co-equalizers (see (\ref{coeqlim})).
\end{enumerate}
\end{proof}

\begin{definition}[Compactly supported analytic index]
The morphism given by the precedent theorem and its corollary is called 
{\it{The compactly supported analytic index of}} $\gr$.
\end{definition}

In the next subsection we are going to slightly modify the indices we constructed, but before that, we give two important properties of these indices. The first one is related to the Bott morphism, so we first describe what the Bott morphism is in our setting. 

Let $\gd:= \gr \times \Rr^2 \rightrightarrows \go \times \Rr^2$ be the product groupoid, where $\Rr^2 
\rightrightarrows \Rr^2$ is the identity groupoid. Let us first recall that the algebra $\ci_c(\Rr^2)$ (with the punctual product) is stable under holomorphic calculus in $C_0(\Rr^2)$, hence the inclusion induces an isomorphism in $K$-theory. In particular, the Bott element can be seen as an element in $K_0(\ci_c(\Rr^2))\approx K_0(C_0(\Rr^2))$. Therefore, we can consider, for each $0\leq k\leq \infty$, the Bott morphism 
$$Bott_k:K_0(\ck)\rightarrow K_0(C_{c}^{k}(\gd)),$$
that is just the product by the Bott element.

Now, from the fact that the product in $K$-theory is natural it follows that the morphism $Bott_k$ passes to the quotient $K_{0}^{h,k}(\gr)$, {\it i.e.}, we get a Bott morphism
$K_{0}^{h,k}(\gr)\stackrel{Bott}{\longrightarrow}K_{0}^{h,k}(\gd)$. Furthermore, by using universal properties we easily extend this morphism to $K_0^F(\gr)\stackrel{Bott}{\longrightarrow}K_{0}^{F}(\gd)$. The following compatibility result will be very useful in te sequel.

\begin{proposition}
The compactly supported index $ind_{a,\gr}^{F}$ is compatible with the Bott morphism,
$\it{i.e.}$, the following diagram is commutative
\[
\xymatrix{
K_0(A^*\gr) \ar[rr]^-{ind_{a,\gr}^{F}} \ar[d]_-{Bott}
& & K_{0}^{F}(\gr) \ar[d]^-{Bott_F} \\
K_0(A^*\gd) \ar[rr]_-{ind_{a,\gd}^{F}} & & K_{0}^{F}(\gd)
}
\]
\end{proposition}

\begin{proof}
It is enough to check that, for each $k\in \Nat$, the following diagram is commutative:
\[
\xymatrix{
K_0(A^*\gr) \ar[rr]^-{ind_{a,\gr}^{h,k}} \ar[d]_{Bott}
& & K_{0}^{h,k}(\gr) \ar[d]^{Bott_k} \\
K_0(A^*\gd) \ar[rr]_-{ind_{a,\gd}^{h,k}} & & K_{0}^{h,k}(\gd)
}
\]
Now, by multiplying again by the Bott element  (seen in $K_0(\ci_c(\Rr^2))$ as above) we have another Bott morphism
$$Bott_T:\src (\gr^T)\longrightarrow \src ((\gd)^T).$$
We use again the fact that the product in $K$-theory is natural and in particular it commutes with the evaluation morphisms to see that we can decompose the last diagram in commutative diagrams in the following way:
\[
\xymatrix{
K^0(A^*\gr) \ar[rrr]^{ind_{a,\gr}^{h,k}} \ar[ddd]_{Bott} & & &
K_{0}^{h,k}(\gr) \ar[ddd]^{Bott_k} \\
 & K_{0}(\src (\gr^T)) \ar[lu]_{e_0} \ar[r]^{e_1^k} \ar[d]_{Bott_T}&
K_0(C_{c}^{k}(\gr)) \ar[ru]^{\pi} \ar[d]^{Bott_k} & \\
& K_{0}(\src ((\gd)^T)) \ar[ld]^{e_0} \ar[r]_{e_1^k} &
 K_0(C_{c}^{k}(\gd)) \ar[rd]^{\pi} & \\
K^0(A^*(\gd)) \ar[rrr]_{ind_{a,\gd}^{h,k}} & & & K_{0}^{h,k}(\gd). \\
}
\]
\end{proof}

The second property is related with the inclusions of open subgroupoids. Let $\gr \rightrightarrows \go$ be a Lie groupoid and
$\hr \rightrightarrows \ho$ be an open subgroupoid. We have the following compatibility result:

\begin{proposition}
The following diagram is commutative:
\[
\xymatrix{
K^0(A^*\hr) \ar[rr]^{ind_{a,\hr}^{F}} \ar[d]
& & K_{0}^{F}(\hr) \ar[d] \\
K^0(A^*\gr) \ar[rr]_{ind_{a,\gr}^{F}} & & K_{0}^{F}(\gr).
}
\]

where the vertical maps are induced from the inclusions by open subgroupoids.
\end{proposition}

\begin{proof}
It is enough to check the proposition for each index of order $k$, $ind_{a}^{h,k}$, for all $k\in \mathbb{N}$.

First, note that $\hr^T \subset \gt$ is an open subset. 
Even more, the algebra inclusion $$\src (\hr^T) \hookrightarrow \src (\gt)$$ commutes with all evaluations.
In particular, the following diagram is commutative
\[
\xymatrix{
K^0(A^*\hr) \ar[d]_j
& K_{0}(\src (\hr^T)) \ar[l]_{e_0} \ar[d]^j \ar[r]^{e_1} & K_{0}(\ci_c(\hr)) \ar[d]^{j}\\
K^0(A^*\gr)  & K_{0}(\src (\gt)) \ar[l]^{e_0}\ar[r]_{e_1}  & K_{0}(\cg)
}
\]
where the morphisms noted by $j$ are induced by the extension by zero outside the open subsets. 
The conclusion is now immediate as we can decompose the diagram in the enouncement of the proposition in the following way
\[
\xymatrix{
K^0(A^*\hr) \ar[rrr]^{ind_{a,\hr}^{h,k}} \ar[ddd]_j & & & 
K_{0}^{h,k}(\hr) \ar[ddd]^j \\
 & K_{0}(\src (\hr^T)) \ar[lu]^{e_0} \ar[r]^{e_1} \ar[d]_j&
K_0(C_{c}^{\infty}(\hr)) \ar[d]^j \ar[ru]^{\pi} & \\
& K_{0}(\src (\gr^T)) \ar[ld]^{e_0} \ar[r]_{e_1} &
K_0(C_{c}^{\infty}(\gr)) \ar[dr]_{\pi} &  \\
K^0(A^*\gr) \ar[rrr]_{ind_{a,\gr}^{h,k}} & & &  K_{0}^{h,k}(\gr), \\
}
\]
\end{proof}



\subsection*{Periodic compactly analytic index}

\begin{definition}[Periodic compactly analytic index] 
As we saw above, we can consider the Bott morphism 
$K_{0}^{F}(\gr)\stackrel{Bott}{\longrightarrow}
K_{0}^{F}(\gd)$. We can take the inductive limit $$\varinjlim_{m} K_{0}^{F}(\gr \times \Rr^{2m})$$ 
induced by $K_{0}^{F}(\gr \times \Rr^{2m})\stackrel{Bott}{\longrightarrow}
K_{0}^{F}(\gr \times \Rr^{2(m+1)})$. We note this group by 
\begin{equation}\label{KBF}
K_{0}^{B}(\gr)=\varinjlim_{m} K_{0}^{F}(\gr \times \Rr^{2m}),
\end{equation}
and we call it the {\it Periodic $K-$theory of }$\gr$.

Let $$ind_{a,\gr}^{B}:K^0(A^*\gr)\rightarrow K_{0}^{B}(\gr)$$ be the morphism given by the composition of $ind_{a,\gr}^{F}:K^0(A^*\gr)\rightarrow
 K_{0}^{F}(\gr)$ followed by $K_{0}^{F}(\gr)
 \stackrel{Bott}{\longrightarrow}K_{0}^{B}(\gr)$. We call this morphism the 
{\it Periodic compactly analytic index of} $\gr$.
\end{definition}

\begin{remark}\label{perbott}
\begin{enumerate}

\item $K_{0}^{B}(\gr)$ satisfy Bott periodicity by construction, {\it i.e.}, 
$$K_{0}^{B}(\gr) \stackrel{Bott}{\longrightarrow}
K_{0}^{B}(\gd)$$ is an isomorphism.

\item The periodic analytic index is also intermediate between $\ci$ and $C^*_r$, {\it i.e.}, it satisfy diagram (\ref{eldiagrama}) too:
\begin{equation}\label{eldiagramaBott}
\xymatrix{
Ell(\gr) \ar[d]_{symb} \ar[r]^{ind}&
K_0(\ci_c(\gr)) \ar[d] & \\ K_0(A^*\gr) \ar[d]_{id} \ar[r]^{ind_{a}^{B}}& K_{0}^{B}(\gr) \ar[d] & \\
K_0(A^*\gr) \ar[r]_{ind_a}& K_0(C_{r}^{*}(\gr)) & 
.}
\end{equation}
where $K_{0}^{B}(\gr)\rightarrow K_0(C_{r}^{*}(\gr))$ 
is induced from
$$K_{0}^{F}(\gr \times \Rr^{2m})\rightarrow 
K_0(C_r^*(\gr \times \Rr^{2m}))\cong K_0(C_{r}^{*}(\gr))$$ 
using the Bott periodicity of $K$-theory for $C^*$-algebras.
\end{enumerate}
\end{remark}

\section{Applications}

\subsection{Longitudinal index theorem}

Let $(M,F)$ be a foliated manifold with holonomy groupoid $\gr \rightrightarrows
M$. In this case the Lie algebroid is given by the integrable subbundle $F$.

In \cite{CS}, Connes-Skandalis define a topological index
$ind_t:K^0(F^*)\rightarrow K_0(C_{r}^{*}(\gr))$ and they show the equality with the 
$C^*$-analytic index of $\gr$. 

We will establish a more primitive longitudinal index theorem. That is, we will see that the equality between the indices takes already place in the group $K_0^B(\gr)$. 

Before stating the longitudinal index theorem, we will need the next proposition (for more details see \cite{Concg} II.5 or \cite{DLN} section 6.1).

\begin{proposition}\label{propDLN}
Let $N\stackrel{p}{\longrightarrow} T$ be a vector bundle over $T$. We consider the Thom groupoid associated to it, {\it i.e.}, 
$\hr:=N\times_T N\rightrightarrows N$, which has Lie algebroid
$N\oplus N^*\cong N\otimes \mathbb{C}$. Then, the following diagram is commutative:
\[
\xymatrix{
K^0(A^*\gr_{\mathscr{T}})\ar[d]_{ind_{a,\gr_{\mathscr{T}}}} \ar[r]^{Thom^{-1}}
 & K^0(T) \ar[d]^{\mathscr{M}} \\
K_0^F(\gr_{\mathscr{T}})& K_0(\ci_c(\gr_{\mathscr{T}})) \ar[l]^{\pi},
}
\]
where $\mathscr{M}$ is the morphism given by the Morita equivalence between $\gr_{\mathscr{T}}$ and $T$. 
In other words, modulo Fourier and Morita, the compactly supported analytic index of $\gr_{\mathscr{T}}$
coincides with the Thom isomorphism's.
\end{proposition}

\begin{proof}
It is known (see for example \cite{DLN} theorem 6.2 or \cite{Concg} II.5) that the following diagram is commutative:
\[
\xymatrix{
K^0(A\hr)\ar[r]^{Thom^{-1}} & K^0(T) \ar[d]^{\mathscr{M}}_{\approx} \\
K^0(A^*\hr)\ar[r]_{ind_{a,\hr}}
\ar[u]^{Fourier}_{\approx} & K_0(C_{r}^{*}(\hr)),
}
\]
where $\mathscr{M}$ is the morphism given by the Morita equivalence between $\gr_{\mathscr{T}}$ and $T$.

It will be then enough to prove that $\ci_c(\hr)\subset C_r^*(\hr)$ is stable under holomorphic calculus. Now, from the classical fact that the algebra of smooth kernel operators $\ci_c(\Rr^{2m} \times \Rr^{2m})$ is stable under holomorphic calculus on the compact operators algebra $\Kom$, we easily get that $\ci_c(\hr)$ is stable under holomorphic calculus on $C^*_r(\hr)$ (locally it reduces to the case of smooth kernel operators, see \cite{trepdo} for instance). In particular $K_0(\ci_c(\hr))\approx K_0(C^*_r(\hr))$ and the proof is complete.
\end{proof}

Let us now define the periodic topological index of a foliation. The definition is analogue to the Connes-Skandalis definition of the $C^*$-topological index.

Before stating the definition, let us recall the following fact: If $T_V$ is an open transversal of a foliated manifold $(V,F_V)$ with holonomy groupoid $\gr_V$, then there a well defined morphism 
\begin{equation}\label{Tmorph}
K^0(T_V)=K_0(\ci_c(T_V))\stackrel{i}{\longrightarrow}K_0(\ci_c(\gr_V)),
\end{equation}
induced by the inclusion of the open subgroupoid resulting from the restriction to the transversal and a suitable Morita équivalence, \cite{CS,BH}.

\begin{definition}\label{it}[Periodic topological index]
Let $g:M\hookrightarrow \Rr^{2m}$ be an embedding, we consider the foliation 
$M\times \Rr^{2m}$ given by the integrable vector bundle $F\times \{0\}$. 
This foliation has $\tilde{\gr}=\gr \times \Rr^{2m}$ as a holonomy groupoid. 
Let $T$ be the normal vector bundle to the foliation in $\Rr^{2m}$, $T_x:=(g_*(F_x))^{\bot}$. Now, the map $h:T\rightarrow M\times \Rr^{2m}$ given by
$(x,\xi) \mapsto (x,g(x)+\xi)$ allows to identify $T$ with an open transversal of 
$(M\times \Rr^{2m},\tilde{F})$, that we still denote by $T$. 
Let $N$ be the normal vector bundle to the inclusion $T\subset M\times \Rr^{2m}$, we can take a neighborhood
$\Omega$ of $T$ in $M\times \Rr^{2m}$ in such a way that
$\hr:= \tilde{\gr}|_{\Omega} \approx N\times_T N$ where
$N \times_T N \rightrightarrows N$ is the pair groupoid over $T$ (we keep the notation from proposition \ref{propDLN}).
This last groupoid has Lie algebroid $A\hr=N\oplus N\approx F\oplus \Rr^{2m}$. We can then consider
the Bott isomorphism $$K^0(F)\stackrel{Bott}{\rightarrow} K^0(A\hr).$$
\noindent
By the {\it Periodic topological index of $\gr$} we mean the morphism
$$ind_{t,\gr}^{B}:K^0(F)\rightarrow K_{0}^{B}(\gr)$$
given by the composition
\[
K^0(F)\stackrel{Bott}{\longrightarrow} K^0(A\hr)
\stackrel{Thom^{-1}}{\longrightarrow}
K^0(T) \stackrel{\iota}{\longrightarrow} K_{0}^{B}(\gr)
\]
where $\iota$ is given by the composition
\begin{equation}\label{itBH}
\hbox{ \footnotesize $
K^0(T) \stackrel{i}{\longrightarrow} K_{0}(\ci_c(\gr \times
\Rr^{2m}))
\longrightarrow K_{0}^{F}(\gr \times
\Rr^{2m}) \longrightarrow K_{0}^{B}(\gr),
$}
\end{equation}
The morphism $$K^{0}(T) \stackrel{i}{\longrightarrow}
K_{0}(\ci_c(\gr \times
\Rr^{2m}))$$ is the one described in (\ref{Tmorph}) above.
\end{definition}

\begin{remark}
If $K_0^B(\gr)\stackrel{i_B}{\longrightarrow} K_0(C_r^*(\gr))$ is the morphism induced by the inclusion, then 
$$i_B\circ ind_{t,\gr}^{B}=ind_{t,\gr},$$
where $ind_{t,\gr}$ is the $C^*$-topological index of Connes-Skandalis.
\end{remark}

Now, we state the theorem in our setting.

\begin{theorem}\label{til}
Let $(M,F)$ be a foliated manifold. With the same notations as above we have that
$$ ind_{a,\gr}^{B}=ind_{t,\gr}^{B}.$$
In particular $ind_{t,\gr}^{B}$ does not depend on the choices made for its definition.
\end{theorem}

\begin{proof}
We use the same notations as in the definition of the periodic topological index.
We are going to show that the following two morphisms coincide
\begin{align}\label{iatrun}
Bott\circ ind_{a,\gr}^{F}:K^0(F)\rightarrow K_{0}^{F}(\gr \times
\Rr^{2m}),
\end{align}
and
\begin{align}\label{ittrun}
K^0(F)\stackrel{Bott}{\longrightarrow} K^0(A\hr)
\stackrel{Thom^{-1}}{\longrightarrow}
K^0(T) \stackrel{\iota}{\longrightarrow} K_{0}^{F}(\gr \times
\Rr^{2m})
\end{align}
where $\iota$ denotes the composition
\[
K^0(T) \stackrel{\mathscr{M}}{\longrightarrow}
K_{0}(\ci_c(\hr))\longrightarrow K_{0}(\ci_c(\gr\times
\Rr^{2m}))
\longrightarrow K_{0}^{F}(\gr \times
\Rr^{2m}),
\]
that is, we will see that the equality of the indices happens before taking the Bott limit.

Now, from the compatibility of the compactly supported index with open subgrupoids and with Bott morphism, seen in last section, we have that
$ind_{a,\hr}^{F}$ and $ind_{a,\gr \times \Rr^{2m}}^{F}$
coincide modulo the inclusions by open subgroupoids, and we have also that 
$ind_{a,\gr \times \Rr^{2m}}^{F}$ and
$ind_{a,\gr }^{F}$ coincide modulo Bott. 
Hence, all the following diagrams are commutative
\[
\xymatrix{
K^0(F) \ar[ddd]_{ind_{a,\gr}^{F}} \ar[rr]^{Bott} \ar[rd]_{Bott}& &
K^0(A\hr)\ar[dd]_{ind_{a,\hr}^{F}} \ar[ld]_{\approx}^i
\ar[r]^{Thom^{-1}} & K^0(T) \ar[dd]^{\mathscr{M}} \\
&K^0(A\gr \times \Rr^{2m}) \ar[rdd]_{ind_{a,\gr \times \Rr^{2m}}^{F}}&
& \\
& & K_{0}^{F}(\hr)\ar[d]^i & K_0(\ci_c(\hr )) \ar[d]^i
\ar[l]^{\pi}\\
K_{0}^{F}(\gr)\ar[rr]_{Bott} & & K_{0}^{F}(\gr \times \Rr^{2m}) & K_0(\ci_c(\gr
\times \Rr^{2m})) \ar[l]^{\pi},
}
\]
where the top right square is commutative from proposition \ref{propDLN}.
To conclude, we have just to remark that in the precedent diagram we have precisely, on one side, the morphism (\ref{iatrun}), left and first below; and on the other side the morphism (\ref{ittrun}).
\end{proof}

\subsubsection{Bounded assembly map}

We recall that the map $$D\mapsto ind_a(D)\in
K_0(C_r^*(\gr))$$ allows to construct an assembly map
\begin{equation}\label{muintro}
\mu:K_{*,\tau}(B\gr)\rightarrow K_0(C_r^*(\gr))
\end{equation}
by putting $\mu(\delta_D)=ind_a(\sigma_D)$. The fact that it is well defined can be deduced (at least in the case of foliations) from the Connes-Skandalis longitudinal index theorem, \cite{CS}.
This morphism was first defined by Baum and Connes \cite{BC} for groups (see \cite{Tu3} for the case of groupoids). 

In our setting, the reinforced longitudinal index theorem allows us to define 
the corresponding assembly map,
\begin{equation}\label{muFintro}
\mu_F:K_{*,\tau}(B\gr)\rightarrow K_{0}^{B}(\gr),
\end{equation}
given as in (\ref{muintro}) but with periodic analytic index: 
$\mu_F(\delta_D)=ind_{a}^{B}(\sigma_D)$. By definition we have a commutative diagram
\begin{equation}\label{musintro}
 \xymatrix{
K_{*,\tau}(B\gr) \ar[r]^-{\mu} \ar[rd]_-{\mu_F} & K_0(C_r^*(\gr)) \\
& K_{0}^{B}(\gr) \ar[u]_i.
}
\end{equation}

\begin{remark}
We recall that the Baum-Connes conjecture establish that $\mu$ is an isomorphism. Its importance stands in the several consequences it would have, \cite{BCH,KL}. 
\end{remark}


\subsection{Pairings with Cyclic cohomology}
 
As we mentioned in the introduction, the main motivation for constructing the compactly supported indices is that, at this level, we can extract numerical information from this $K$-theory elements. There is 
indeed a pairing between $K-$theory and Cyclic cohomology (see (\ref{accouplement})), and an important problem in non commutative geometry is to give explicit (topological) formulae for this pairings.

Now, for fixed $\tau \in HP^*(\cg)$, we can expect an easy (topological) calculation only if the map 
$D\mapsto \langle D \, , \tau \rangle$ factors through the symbol class of $D$, $[\sigma(D)]\in K^0(A^*\gr)$. We are going to solve the factorization problem for bounded cocycles. This restriction is not at all restrictive. The bounded cocycles are periodic cyclic cocycles whose formulas only use a finite number of derivates. This is of course the case of group (and group action) cocycles, the transverse fundamental class, Godbillon-Vey and all the cocycles coming from $H^*(B\gr)$ (as we see in the last section).  In fact, as we will see below, we are going to completely solve the problem for {\'e}tale groupoids. Before stating the factorization theorem, we give the precise definition of bounded cocycles. 

\begin{definition}\label{defBCC}
A multilinear map $\tau:\underbrace{\cg \times \cdots \times \cg}_{q+1-times}
\rightarrow \mathbb{C}$ is bounded if it extends to a continuous multilinear map 
$\underbrace{\ck \times \cdots \times \ck}_{q+1-times}\stackrel{\tau_k}{\longrightarrow}\mathbb{C}$, for some $k\in \mathbb{N}$. We can consider a sub-bicomplex of the Periodic bicomplex $(C^{n,m}(\cg),b,B)$ consisting in bounded multilinear maps. We note this bicomplex by $(C_{B}^{n,m}(\cg),b,B)$ and a cocycle for it will be called a bounded continuous cyclic cocycle.
\end{definition}

In the following proposition we prove that for the pairing with a bounded cocycle extends to our groups 
$K_0^B(\gr)$.

\begin{proposition}\label{extfin}
Let $\tau$ be a bounded cyclic cocycle. Then the pairing map 
$K_0(\cg)\stackrel{\langle \, ,\tau \rangle}{\longrightarrow}$ extends to $K_0^B(\gr)$, {\it i.e.}, we have a commutative diagram of the following type:
\begin{equation}\label{pairF}
\xymatrix{
K_0(\cg) \ar[r]^-{<,\tau>} \ar[d]_-{\iota} & \mathbf{C} \\
K_0^B(\gr) \ar[ru]_-{\tau_B} &
}
\end{equation}
\end{proposition}

\begin{proof}
We are going to see first that the pairing extends to $K_0^F(\gr)$.

Let $\tau_k: \underbrace{\ck \times \cdots \times \ck}_{2n+1-times}\rightarrow \mathbb{C}$ be a continuous extension of $\tau$, for some $k\in \mathbb{N}$. It is immediate by definition that $K_0(\cg)\stackrel{\langle \, ,\tau \rangle}{\longrightarrow}\mathbb{C}$ extends to $K_0(\ck)$, more explicitly, the following diagram is commutative
\[
\xymatrix{
K_0(\cg) \ar[r]^-{<,\tau>} \ar[d]_-{(\iota_k)_*} & \mathbf{C} \\
K_0(\ck) \ar[ru]_-{<,\tau_k>} &
}
\]
We are now going to see that the pairing passes to $K_{0}^{h,k}(\gr)$. For that we just have to check that this pairing preserves the relation $\sim_h$ over $K_0(\ck)$. This can be done by adapting an argument already used by Connes \cite{Concdg} and Goodwillie \cite{Good}: let $e\in \ci([0,1],\ck)\oplus \mathbb{C}$ be an idempotent. It defines a smooth family of idempotents $e_t$ in $\widetilde{\ck}$. 
We set $a_t:=\frac{de_t}{dt}(2e_t-1)$. 
Hence, a simple calculation shows 
$$\frac{d}{dt} \langle \tau ,e_t \rangle= \sum_{i=0}^{2n}\tau(e_t,...,[a_t,e_t],...,e_t)=:L_{a_t}\tau(e_t,...,e_t).$$
Now, {\it{the Lie derivates}} $L_{x_t}$ act trivially on $HP^0(\ck)$ (see \cite{Concdg,Good}), then $ \langle \tau ,e_t \rangle$ is constant in $t$. In particular, $ \langle \tau ,e_0 \rangle= \langle \tau ,e_1 \rangle$. It follows immediately that $\langle \, ,\tau \rangle$ extends to 
$K_{0}^{F}(\gr):=\varinjlim_kK_{0}^{h,k}(\gr)$.

Finally, the Periodic Cyclic cohomology $HP^*(\cg)$ satisfies Bott periodicity. Hence, the extension from 
$K_{0}^{F}(\gr)$ to $K_{0}^{B}(\gr):=\varprojlim_{m}K_{0}^{F}(\gr \times \Rr^{2m})$ is now immediate.
\end{proof}

\begin{remark}\label{remextfin}
\quad
\begin{enumerate}
\item The extension of diagram (\ref{pairF}) is very explicit: let $\tau \in C^{n,m}_{F}(\gr)$ and $x=(x_1,x_2,...)\in
K_0^F(\gr)=\varprojlim_kK_{0}^{h,k}(\gr)$. Let $\tau_k$ be an extension
of $\tau$. Then 
$$\tau_F(x)= \langle \tau_k ,x_k \rangle.$$
\item If we note by $HP_{B}^{*}(\ci_c(\gr))$ the cohomology of the bicomplex $(C_{B}^{n,m}(\cg),b,B)$ described in the definition \ref{defBCC}, then we have a pairing
\begin{equation}\label{pairBottF}
HP_{B}^{*}(\ci_c(\gr))\times K_{0}^{B}(\gr)\longrightarrow \mathbb{C},
\end{equation}
induced from the extension of the proposition above. More explicitly, if $\tau$ is a bounded cocycle and $x\in K^F_0(\gr)$, then $$\langle x,\tau \rangle=\langle Bott(x),Bott(\tau)\rangle.$$
\end{enumerate}
\end{remark}

\begin{theorem}[Factorization theorem]\label{extfinteo}
\quad
\begin{itemize}
\item[{\it (a)}] Let $\tau $ be a bounded continuous cyclic cocycle. Then $\tau$ defines a morphism 
$\varphi_{\tau}:K_0^B(\gr)\rightarrow \mathbb{C}$ such that the following diagram is commutative
\[
\xymatrix{
Ell(\gr) \ar[r]^-{ind} \ar[d]_-{symb.} & K_0(\cg) \ar[d] \ar[rr]^-{\langle \_ , \tau \rangle} & & \mathbb{C} \\
K^0(A^*\gr) \ar[r]_-{ind_{a}^{B}}& K_{0}^{B}(\gr) \ar[urr]_-{\varphi_{\tau}} &  &.
}
\]
In particular, for a 
$\gr$-pseudodifferential elliptic operator $D$, we have the following formula
\begin{equation}
\varphi_{\tau}\circ ind_{a}^{B}([\sigma(D)])=<ind\, D,\tau>.
\end{equation}
\item[{\it (b)}] Let $\gr \rightarrow \go$ be an {\'e}tale groupoid. Then we have the result of precedent paragraph for every $\tau \in HP^*(\cg)$. In particular, the map 
$$Ell(\gr)\stackrel{ind}{\longrightarrow}K_0(\cg) \stackrel{\langle \_ ,\tau \rangle}{\longrightarrow} \mathbb{C}$$ always factors through $Ell(\gr)\rightarrow K^0(A^*\gr)$.
\end{itemize}
\end{theorem}

\begin{proof}
\quad
\begin{itemize}
\item[{\it (a)}] It is immediate from proposition \ref{extfin}.
\item[{\it (b)}] Thanks to the works of Burghelea, Brylinski-Nistor and Crainic (\cite{Bur,BN,Cra}), we known a very explicit description of the Periodic cyclic cohomology for {\'e}tale groupoids. For instance, we have a decomposition of the following kind (see for example \cite{Cra} theorems 4.1.2. and 4.2.5)
\begin{equation}\label{Crainic}
HP^*(\cg)=\Pi_{\ops}H_{\tau}^{*+r}(B\Nb_{\ops}),
\end{equation}
where $\Nb_{\ops}$ is an {\'e}tale groupoid associated to $\ops$ (the normalizer of 
$\ops$, see 3.4.8 in ref.cit.). For instance, when $\ops=\go$, $\Nb_{\ops}=\gr$.

Now, all the cyclic cocycles coming from the cohomology of the classifying space are bounded. Indeed, we know that each factor of 
$HP^*(\cg)$ in the decomposition (\ref{Crainic}) consists of bounded cyclic cocycles (see the last section of this work). 
Now, the pairing 
$$HP^*(\cg)\times K_0(\cg)\longrightarrow \mathbb{C}$$
is well defined. In particular, the restriction to $HP^*(\cg)|_{\ops}$ vanishes for almost every $\ops$. The conclusion is now immediate from proposition \ref{extfin}.
\end{itemize}
\end{proof}


\subsubsection{Geometric Corollaries}

The factorization problem we just have met is deeply related with the Novikov conjecture. Indeed, if the map
$$K_0(\cg)\stackrel{\langle \_ ,\tau \rangle}{\longrightarrow} \mathbb{C}$$
where $\tau \in HP^*(\cg)$, extends to $K_0(C^*_r(\gr))$, then the factorization through the principal symbol class is immediate. However, as the following example shows, it is far from being a trivial problem. 
    
\begin{ejemplo}\label{intronov}\cite{CMnov,Concg}
Let $\Gamma$ be a discrete group acting properly and freely on a smooth manifold $\tilde{M}$ with compact quotient
$\tilde{M}/\Gamma:=M$. Let $\gr\rightrightarrows \go=M$ be the Lie groupoid quotient of $\tilde{M}\times
\tilde{M}$ by the diagonal action of $\Gamma$. 

Let $c \in H^{*}(\Gamma):=H^*(B\Gamma)$. Connes-Moscovici showed in \cite{CMnov} that the higher Novikov signature, $Sign_c(M)$, 
can be obtained with the paring of the signature operator $D_{sign}$
and a cyclic cocycle $\tau_c$ associated to $c$: 
\begin{equation}\label{signind}
\langle \tau_c , ind\, D_{sgn}\rangle=Sign_c(M,\psi). 
\end{equation}
The Novikov conjecture states that these higher signatures are oriented homotopy invariants of
$M$. Hence, if
$ind\, D_{sign}\in K_0(\cg)$ is a homotopy invariant of
$(M,\psi)$ then the Novikov conjecture would follow. We only known that 
$j(ind\, D_{sign})\in K_0(C^*_r(\gr))$ is a homotopy invariant. But then we have to extend the action of $\tau_c$ to $K_0(C^*_r(\gr))$. Connes-Moscovici show that this action extends for Hyperbolic groups.
\end{ejemplo}

{\bf{Question }} Is $ind_{a}^{B}(D)$ a homotopy invariant? An affirmative answer to this question would imply the Novikov conjecture because the pairing (\ref{signind}) extends to $K_0^B(\gr)$. In fact, since we know that $ind_a(D)\in K_0(C^*_r(\gr))$ is a homotopy invariant, another way to establish the Novikov conjecture would be to prove the injectivity of the map $K_0^B(\gr)\stackrel{i}{\longrightarrow}K_0(C^*_r(\gr))$.
\bigskip
 
In \cite{Cotfc}, theorem 8.1, Connes solves the extension problem for some kind of cyclic cocycles over the holonomy groupoid of a foliation, and he gives a topological formula for the pairing. 
A main step in his proof is the Connes-Skandalis longitudinal index theorem.

Using our longitudinal index theorem we obtain, as a corollary of Connes theorem, the  analog result in our setting. In this case we do not have to deal with the extension problem, and so the result applies to all classes $c\in H^*(B\gr)$.  

\begin{corollary}
Let $(V,F)$ be a foliated manifold (non necessarily compact) transversally oriented. 
Let $\gr$ be its holonomy groupoid. 
For any $c\in H^*(B\gr)$ there is an additive map 
$$\varphi_c:K_{0}^{B}(\gr)\rightarrow \mathbb{C}$$ such that
\begin{equation}
\varphi_c(\mu_F(x))
=\langle ch_{\tau}(x),c\rangle \, \, \forall x\in K_{0}^{B}(\gr).
\end{equation}
\end{corollary}

\subsubsection*{Haefliger Cohomology}

The reinforced longitudinal theorem can also be used for re-establish the index formulas in Haefliger cohomology found by Benameur-Heitsch in \cite{BH}. 

Benameur-Heitsch start by defining an algebraic Chern character
\begin{equation}
ch_a:K_0(\cg)\rightarrow H_{c,bas}^{*}(M/F),
\end{equation}
where $H_{c,bas}^{*}(M/F)$ is the Haefliger cohomology (see \cite{Haef} or \cite{CrMoercr}).
This character is compatible with the
``shriek maps" in $K$-theory and in Haefliger cohomology. 

The main result in \cite{BH} requires the construction of a map 
$ch_a(ind_t):K^0(T^*F)\rightarrow H_{c,bas}^{*}(M/F)$ that is morally the topological index followed by the Chern character. Then, they prove the following formula (theorem 5.11 ref.cit.): 

For any $u\in K^0(T^*F)$,
\begin{equation}\label{chatintro}
ch_a(ind_t)(u)=(-1)^p \int_{F}\pi_F!(ch(u))Td(F\otimes \mathbb{C})
\in H_{c,bas}^{*}(M/F),
\end{equation}
where $\pi_F!:H^*_c(F, \Rr)\rightarrow H^*(M,\Rr)$ is the fiberwise integration, and $ch:K^0(F)\rightarrow H^*_c(F)$ is the usual Chern character. 

Using the same kind of arguments as in proposition \ref{extfin}, one can show that the character $ch_a$ extends to the group $K_{0}^{B}(\gr)$:
\[
\xymatrix{
K_0(\cg) \ar[r] \ar[d]_-{ch_a} & K_{0}^{B}(\gr) \ar[dl]^-{ch_a} \\
H_{c,bas}^{*}(M/F) & .
}
\]

The next result is an immediate consequence of the reinforced longitudinal theorem and the Benameur-Heitsch formula.

\begin{corollary}\label{corBHintro}
For any longitudinal elliptic pseudodifferential operator $D$, we have that
\begin{equation}
ch_a(ind\, D)=ch_a(ind_t)[\sigma_D]
\in H_{c,bas}^{*}(M/F).
\end{equation}
\end{corollary}

The fact that the actions of a large class of cyclic cocycles extend naturally to the periodic group 
$K_{0}^{B}(\gr)$ allows to think that formulas as those of Benameur-Heitsch 
could be developed in cohomology spaces more complex than $H_{c,bas}^{*}(M/F)$.

\subsubsection{Bounded Cyclic cohomology}
Let $HP^*_B(\cg)$ be the cohomology of the bicomplex $(C_{B}^{n,m}(\cg),b,B)$ 
(see definition \ref{defBCC} for notations),  we call it the {\it "Bounded (Periodic) Cyclic Cohomology of $\gr$"}.

\subsection*{The case of {\'E}tale groupoids}\label{etaleapp}

Let us consider the classifying space, $B\gr$, and its twisted cohomology, $H_{\tau}^{*}(B\gr)$. 

In \cite{Cotfc}, Connes constructed
a group morphism (defined at the level of cocycles)
\begin{equation}\label{tordue}
\phi: H_{\tau}^{*}(B\gr) \rightarrow HP^*(\cg),
\end{equation}
that is in fact an inclusion as a direct factor (see also \cite{Concg}, $III.2.\delta$).
 
By examinining Connes' construction, it follows easily that the cocycles one gets are all bounded. In other words, one constructs a morphism (\cite{Ca3})
\begin{equation}\label{tordueF}
\phi_B: H_{\tau}^{*}(B\gr) \rightarrow HP^*_B(\cg),
\end{equation}
fitting the diagram
\begin{equation}\label{torduediag}
\xymatrix{
H_{\tau}^{*}(B\gr) \ar[r]^-{\phi} \ar[rd]_-{\phi_B} & HP^*(\cg)\\
& HP^*_B(\cg) \ar[u]_-{\iota}.
}
\end{equation}
This means that the image of the morphism in (\ref{tordue}) consists only on Bounded Cyclic cocycles.

\bibliographystyle{amsalpha}
\bibliography{bibliographie}

\end{document}